\documentclass[12pt]{amsart}

\usepackage{amssymb,latexsym}

\usepackage{graphicx,epsfig}
\usepackage[dvips]{color}

\def\cal{\mathcal}

\def\Bbb{\mathbb}

\def\rank{\text{\rm rank\,}}

\def\ord{\text{\rm ord\,}}

\def\pad{\phi^a}

\def\ra{\rangle}
\def\laa{\langle}

\def \supp {\text{\rm supp\,}}
\def\dist{\text{\rm dist\,}}

\def\hl{{h_{\rm \,lin}}}
\def\A{{\cal A}}
\def\D{{\cal D}}
\def\F{{\cal F}}
\def\M{{\cal M}}
\def\N{{\cal N}}
\def\cS{{\cal S}}
\def\T{{\cal T}}
\def\bC{{\Bbb C}}

\def\NN{{\Bbb N}}
\def\bN{{\Bbb N}}
\def\bR{{\Bbb R}}
\def\RR{{\Bbb R}}

\def\al{{\alpha}}
\def\be{{\beta}}
\def\ga{{\gamma}}
\def\la{{\lambda}}
\def\La{{\Lambda}}
\def\om{{\omega}}

\def\x{(x_1,x_2)}
\def\y{(y_1,y_2)}
\def\pa{{\partial}}

\def\ve{{\varepsilon}}
\def\si{{\sigma}}

\def\de{{\delta}}

\def\Om{{\Omega}}
\def\ka{{\kappa}}

\def\th{{\theta}}

\def\pr{\text{\rm pr\,}}

\def\bpm{\begin{pmatrix}}
\def\epm{\end{pmatrix}}

\def\noi{\noindent}
\def\bee{\begin{enumerate}}
\def\ee{\end{enumerate}}

\def\qed{\smallskip\hfill Q.E.D.\medskip}

\newtheorem{thm}{Theorem}[section]
\newtheorem{namedthm}[thm]{Theorem}
\newtheorem{prop}[thm]{Proposition}
\newtheorem{proposition}[thm]{Proposition}
\newtheorem{cor}[thm]{Corollary}
\newtheorem{lemma}[thm]{Lemma}
\newtheorem{remark}[thm]{Remark}
\newtheorem{remarks}[thm]{Remarks}

\newtheorem{assumption}[thm]{Assumption}
\newtheorem{example}[thm]{Example}

\begin{document}

\title[harmonic analysis and hypersurfaces]{
Problems of Harmonic Analysis related  to finite type hypersurfaces in $\bR^3,$ and Newton polyhedra}

%\author[I. A. Ikromov]{Isroil A. Ikromov}
%\address{Department of Mathematics, Samarkand State University,
%University Boulevard 15, 140104, Samarkand, Uzbekistan}
%\email{{\tt ikromov1@rambler.ru}}

\author[D. M\"uller]{Detlef M\"uller}
\address{Mathematisches Seminar, C.A.-Universit\"at Kiel,
Ludewig-Meyn-Stra\ss{}e 4, D-24098 Kiel, Germany}
\email{{\tt mueller@math.uni-kiel.de}}
\urladdr{{http://analysis.math.uni-kiel.de/mueller/}}

\thanks{2010 {\em Mathematical Subject Classification.}
42B25, 42B37, 42B99}
\thanks{{\em Key words and phrases.}
Hypersurface, oscillatory integral, oscillation index, maximal function, sub-level estimate,  Fourier restriction, Newton diagram}

\begin{abstract}
This article is intended to give an overview  on  a collection of results which have been obtained  jointly with I.I. Ikromov, and in parts also with M. Kempe, and at the same time to give a kind of guided tour through the rather comprehensive proofs of the major results that I shall address.   All of our work is highly influenced by the pioneering ideas developed by E.M. Stein.
\smallskip

 Consider a smooth hypersurface $S$  in $\RR^3$ of finite type, and let  $d\mu=\rho d\si$ be a surface carried measure with smooth, compactly supported density $\rho\ge 0$ with respect to the surface measure $d\si.$ The problems that I shall address 
 are the following ones:
 
A. Find, if  possible, optimal  uniform decay estimates for the Fourier transform of the surface carried measure $d\mu.$

B. If we denote by $A_t$ the  averaging operator
$
A_tf(x):=\int_{S} f(x-ty) d\mu(y),
$  determine for which exponents $p$ the associated maximal operator 
$$
\M f(x):=\sup_{t>0}|A_tf(x)|
$$
is bounded on $L^p(\RR^3).$ 

C. Determine the range of  exponents $p$ for which a Fourier restriction estimate 
$$
\Big(\int_S|\hat f(x)|^2\, d\mu(x)\Big)^{1/2}\le C \|f\|_{L^p(\RR^3)}
$$
holds true. The first problem is classical, and the other ones originate from seminal contributions  by E.M. Stein.   These problems are somewhat interrelated, and in particular it is well-known that the answers to B and C are strongly influenced by the answer to A.  Moreover, due to the insights by V.I. Arnol'd and his school, in particular through work by A.N. Varchenko and V.N. Karpushkin, we know that the answer to problem A can be given in terms of Newton polyhedra associated to the given  hypersurface $S,$ at least for analytic hypersurfaces. That Newton polyhedra would play an important role also  in other problems of harmonic analysis  has been stressed later in pioneering work by D.H. Phong and E.M. Stein on oscillatory integral operators. 

In this article, I shall outline how (almost)  complete answers to  the above questions can be given again  in terms of associated Newton polyhedra (for question B, at least if $p>2$).

\end{abstract}

\maketitle

%\vfill\newpage

\tableofcontents

\thispagestyle{empty}

\setcounter{equation}{0}
\section{Introduction}\label{intro}

Let $S$ be a smooth, finite type  hypersurface in $\RR^3$ with Riemannian 
surface measure $d\si,$ and consider the compactly supported measure 
$d\mu:=\rho d\si$ on $S,$ where $0\le\rho\in C_0^\infty(S).$

The problems on which I shall essentially focus  are the following ones:
\medskip

A.  Find, if possible, optimal uniform decay estimates for the Fourier transform of the surface carried measure $d\mu.$
\medskip

B.  If we denote by $A_t$ the  averaging operator
$
A_tf(x):=\int_{S} f(x-ty) d\mu(y),
$  determine for which exponents $p$ the associated maximal operator 
$$
\M f(x):=\sup_{t>0}|A_tf(x)|
$$
is bounded on $L^p(\RR^3).$ For instance, if $S$ is a Euclidean sphere centered at the origin, then $\M$ is the spherical maximal operator studied first by Stein \cite{stein-sphere}.
\medskip

C.  Determine the range of  exponents $p$ for which a Fourier restriction estimate 
$$
\Big(\int_S|\hat f(x)|^2\, d\mu(x)\Big)^{1/2}\le C \|f\|_{L^p(\RR^3)}
$$
holds true.
I shall explain how these questions can be answered in an (almost)  complete way in terms of Newton polyhedra associated to the given surface $S.$  All these results are based on  joint work with Ikromov, and in parts also with Kempe  \cite{IM-ada}, \cite{IM-uniform},  \cite{ikm}, \cite{IKM-max} \cite{IM-rest1},\cite{IM-rest2}.

\medskip
Problem A is a classical question about estimates for oscillatory integrals, and there exists a huge body of results on it, in particular for convex hypersurfaces (some references, also to higher dimensional results,  will follow). 

\medskip
The other two problems had first been formulated by Stein: the study of maximal averages along  hypersurfaces has been initiated in  Stein's work on the spherical maximal function \cite{stein-sphere}, and also the idea of Fourier restriction  goes back to him (cf. \S \ref{FREST}). Their great importance for the study of partial differential equations became clear through Strichartz' article \cite{strichartz}, and in the PDE-literature their dual versions are often called Strichartz estimates.

Of course, these and related  problems, such as the question of optimal sub-level estimates and integrability indices,  make perfect sense also in higher dimensions (and possibly higher co-dimension), but only partial answers are known to most of these problems in this case. The reason for this is that strong information on the  resolution of singularities is available for analytic functions of two variables, for instance by means of Puiseux series expansions of roots, whereas the situation for multivariate functions of more than two variables  is  substantially more complex. Nevertheless, there has been a lot of progress on the question as to how to construct  more elementary and ''concrete'' resolutions of singularities for real analytic multivariate functions, giving more detailed information than what Hironaka's  celebrated  theorem \cite{hironaka} on the resolution of singularities would yield, for instance in work by Bierstone and P.D. Milman \cite{bierstone-milman1}, \cite{bierstone-milman2},  Sussman \cite{sussman},  Parusi\'nski \cite{parusinski1}, \cite{parusinski2},  Greenblatt \cite{greenblatt-resol1}, \cite{greenblatt-resol2}, Collins, Greenleaf and Pramanik \cite{cgp}, among others.  These techniques have already led   to a very good understanding  of,  for instance, the sublevel estimation problem for slices of the surfaces in direction to the Gaussian normal,  and the related determination of critical integrability indices, in independent work by Greenblatt \cite{greenblatt-resol1}, \cite{greenblatt-resol2}, and also Collins, Greenleaf and Pramanik \cite{cgp}, by rather different methods. 
Yet another approach has been developed   by Denef, Nicaise and Sargos  \cite{dns} in order to estimate oscillatory integrals  in higher dimensions.

These  encouraging results give rise to the hope that eventually, substantial progress might also be possible in higher dimensions for the deeper problems A - C, but at present, only partial answers are available, for restricted classes of  surfaces (see, e.g., \cite{greenblatt1}, \cite{greenblatt-max2}).    What makes these problems indeed a lot harder than the sub-level set problem  is that it will surely require estimates of oscillatory integrals with phase functions depending also on small parameters, not just a fixed phase function. This leads to the fundamental problem of stability of the  estimates in Problem A under small perturbation; I shall briefly come back to this later.

\medskip
Before returning to the two-dimensional case, let me address  some results  on the three
problems  A - C  which have been obtained  for particular classes of hypersurfaces also  in higher dimensions. 

The case of convex hypersurface of finite line type has been studied quite intensively, for instance in work by  Randol \cite{randol},  Svensson \cite{svensson},  Schulz \cite{schulz}, Bruna, Nagel and Wainger \cite{bruna-n-w} and Cowling, Disney, Mauceri and myself  \cite{cowling-d-m} concerning problem A,  and for instance by  Nagel, Seeger and Wainger \cite{nagel-seeger-wainger} and Iosevich, Sawyer and Seeger \cite{io-sa-seeger} on problem B. 
Various particular classes of non-convex hypersurfaces have been examined too, for instance in work by Iosevich and Sawyer \cite{iosevich-sawyer1},  Cowling and Mauceri \cite{cowling-mauceri2}, \cite{cowling-mauceri1},  and by Sogge and Stein \cite{sogge-stein}, where the method of damping (with powers of the Gaussian curvature) had been introduced in order to derive partial answers to problem B for very general hypersurfaces in $\bR^n.$ 

\medskip
However, from now on I shall  concentrate on the two-dimensional case, and give  only  occasionally a few references to work in higher dimensions.  Since all three problems can be seen as ''classical'' by now, there is an abundance of literature associated to them, so that it would seem impossible to give due credit to everyone who has made contributions, and I apologize in advance  to everyone whose work is not mentioned, due to lack of space or my personal ignorance. Many further  references can be found in the cited articles, and for a more detailed account of the state of the theory and its historical development until roughly a decade ago, I refer, for instance,  to  Stein's monograph \cite{stein-book}.

\bigskip
It is obvious that all three problems A - C can be localized to sufficiently small neighborhoods of given points  $x^0$ on $S.$  Observe also that the problems A and B are invariant under translations and rotations of the ambient space, so that we may replace the surface $S$ by any suitable image under a Euclidean motion of $\RR^3.$ 

We may thus assume that $x^0=(0,0,0),$ and that $S$ is the graph 
\begin{equation}\label{graph1}
S=\{(x_1,x_2, \phi\x): \x\in \Om \},
\end{equation}
of a smooth function $\phi$ defined on a  sufficiently small neighborhood $\Om$ of the origin,  such that 
\begin{equation}\label{null}
\phi(0,0)=0,\, \nabla \phi(0,0)=0.
\end{equation}

The situation is  different for problem B, since dilations do not commute with translations, so that we are only allowed to work with linear transformations of the ambient space.  We shall therefore study the maximal operator $\M$ under the following transversality assumption on $S.$

\begin{assumption}\label{s1.1}
The affine tangent plane $x+T_xS$ to $S$ through $x$ does not pass through the origin in $\RR^3$  for every $x\in S.$ Equivalently, $x\notin T_xS$ for every $x\in S,$ so that $0\notin S,$ and $x$ is transversal to $S$ for every point $x\in S.$ 
\end{assumption}

Notice that this assumption allows us to find a linear change of coordinates in $\RR^3$ so that in the new coordinates
 $S$ can locally be represented as the shifted graph of a function $\phi$ as before, more precisely,
 \begin{equation}\label{graph2}
S=\{(x_1,x_2, 1+\phi\x): \x\in \Om \},
\end{equation}
where $\phi$ satisfies again \eqref{null}.

Our transversality assumption is natural in this context. Indeed,  various examples show that if it is not satisfied, then the behavior of the corresponding maximal function may change drastically.

\smallskip
Observe also that if  $\phi$ is flat, i.e., if all derivatives of $\phi$ vanish at the origin, and if $\rho(x^0)>0,$  then it is well-known and easy to see that the maximal operator $\M$ is $L^p$-bounded if and only if $p=\infty,$ so that this case is of no interest.  In a similar way, also the problems A and C will become of quite a different nature. We shall therefore  assume that $\phi$ is non-flat, i.e., of finite type. Correspondingly, we shall always assume that the hypersurface $S$ is of finite type,  in the sense that every tangent plane has finite order of contact. Recall that in the study a convex hypersurface, a standard assumption used by many authors is that the surface is of finite line type, which means that every  tangent line has finite order of contact, which is stronger than our assumption.

\bigskip
The article will be organized as follows:
\medskip

In the next paragraph, I shall briefly review some basic  concepts concerning  Newton polyhedra, adaptedness of coordinates and the notion of height, which have been introduced by Arnol'd (cf. \cite{arnold}, \cite{agv}) and his school, most notably Varchenko \cite{Va}. These concepts had originally been studied for real analytic functions $\phi,$ but as shown in \cite{IM-ada}, can be extended to smooth, finite type functions $\phi.$  Problems A  to C will then be discussed subsequently in $\S 3 - \S 5.$ 

Since many of our  proofs are quite  elaborate, my main goal will be to give a kind of guided tour through the basic structure of the proofs,  hoping that this might also be  helpful for a  everyone interested in studying some of the proofs in more detail.

 \setcounter{equation}{0}
\section{Newton polyhedra, and adapted coordinates}\label{newton}

 Let me first recall
some basic notions from  \cite{Va},\cite{IM-ada}. If $\phi$ is given as before, consider the associated Taylor series 
$$\phi(x_1,x_2)\sim\sum_{\al_1,\al_2=0}^\infty c_{\al_1,\al_2} x_1^{\al_1} x_2^{\al_2}$$
of $\phi$ centered at  the origin.
The set
$$\T(\phi):=\{(\al_1,\al_2)\in\bN^2: c_{\al_1,\al_2}=\frac 1{\al_1!\al_2!}\partial_{ 1}^{\al_1}\partial_{ 2}^{\al_2}­ \phi(0,0)\ne 0\}
$$
will be called the {\it Taylor support } of $\phi$ at $(0,0).$  We shall always assume that
$$\T(\phi)\ne \emptyset,$$
i.e., that the function $\phi$ is of finite type at the origin. The
{\it Newton polyhedron} $\N(\phi)$ of $\phi$ at the origin is
defined to be the convex hull of the union of all the quadrants
$(\al_1,\al_2)+\bR^2_+$ in $\bR^2,$ with $(\al_1,\al_2)\in\T(\phi).$  
The associated {\it Newton diagram}  $\N_d(\phi)$ in the sense of Varchenko
\cite{Va}  is the union of all compact faces  of the Newton
polyhedron; here, by a {\it face,} we shall  mean an edge or a
vertex.

%%%%%%%%%%%%%%%%%%%%%%
%%%%%%%%%%%%%%%%%%%%%%
%%%%%%% Figure 1 %%%%%%%%%
%%%%%%%%%%%%%%%%%%%%%%
%%%%%%%%%%%%%%%%%%%%%%
\begin{figure}
\centering
\scalebox{0.35}{\input{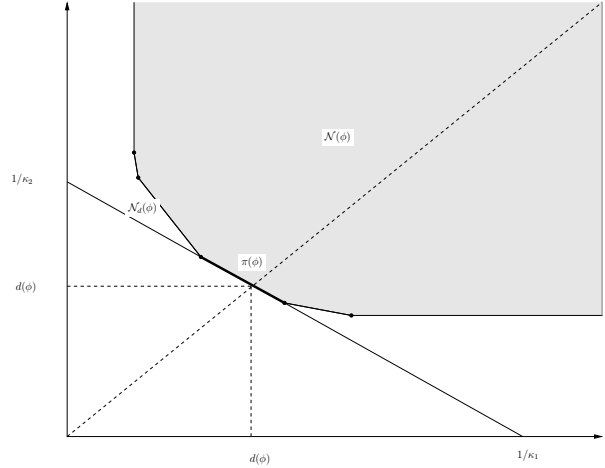}}
\caption{Newton polyhedron}{\label{fig1}}
\end{figure}

We shall use coordinates $(t_1,t_2)$ for points in the plane
containing the Newton polyhedron, in order to distinguish this plane
from the $(x_1,x_2)$ - plane.

The {\it Newton distance}, or shorter {\it distance} $d=d(\phi)$
between the Newton polyhedron and the origin in the sense of
Varchenko is given by the coordinate $d$ of the point $(d,d)$ at
which the bi-sectrix   $t_1=t_2$ intersects the boundary of the
Newton polyhedron.

The {\it principal face} $\pi(\phi)$  of the Newton polyhedron of
$\phi$  is the face of minimal dimension  containing the point
$(d,d)$. Deviating from the notation in \cite{Va}, we shall call the
series
$$\phi_\pr(x_1,x_2):=\sum_{(\al_1,\al_2)\in \pi(\phi)}c_{\al_1,\al_2} x_1^{\al_1} x_2^{\al_2}
$$
the {\it principal part} of $\phi.$ In case that $\pi(\phi)$ is
compact,  $\phi_\pr$ is a mixed homogeneous polynomial; otherwise,
we shall  consider $\phi_\pr$ as a formal power series.

Note that the distance between the Newton polyhedron and the origin
depends on the chosen local coordinate system in which $\phi$ is
expressed.  By a  {\it local  coordinate system at the origin} I 
shall mean a smooth   coordinate system defined near the origin
which preserves $0.$ The {\it height } of the  smooth function
$\phi$ is defined by
$$h(\phi):=\sup\{d_y\},$$
 where the
supremum  is taken over all local  coordinate systems $y=(y_1,y_2)$ at the origin, and where $d_y$
is the distance between the Newton polyhedron and the origin in the
coordinates  $y$.

A given   coordinate system $x$ is said to be
 {\it adapted} to $\phi$ if $h(\phi)=d_x.$

In \cite{IM-ada} we proved that one  can always find an adapted
local  coordinate system in two dimensions, thus  generalizing  the
fundamental work by Varchenko  \cite{Va} who worked in the   setting
of  real-analytic functions $\phi$ (see also \cite{PSS}).

Recall also that if the principal face of the Newton polyhedron $\N(\phi)$ is a compact edge, then it lies on a unique ``principal line'' 
$$
 L:=\{(t_1,t_2)\in \RR^2:\ka_1t_1+\ka_2 t_2=1\}, 
$$
with $\ka_1,\ka_2>0.$ By permuting the  coordinates $x_1$ and $x_2,$ if necessary, we shall always assume that $\ka_1\le\ka_2.$ The  weight $\ka=(\ka_1,\ka_2)$ will be called the {\it principal weight} associated to $\phi.$  It induces dilations $\de_r\x:=(r^{\ka_1}x_1,r^{\ka_2} x_2),\ r>0,$ on $\RR^2,$ so that the principal part $\phi_\pr$ of $\phi$ is $\ka$- homogeneous of degree one with respect to these dilations, i.e.,  $\phi_\pr(\de_r\x)=r\phi_\pr\x$ for every $r>0,$ and 
\begin{equation}\label{dnew}
d=\frac 1{\ka_1+\ka_2}=\frac1{|\ka|}.
\end{equation}
It can then easily be shown (cf. Proposition 2.2 in \cite{IM-ada}) that $\phi_\pr$ can be factorized as 
\begin{equation}\label{2.4}
\phi_\pr\x=c x_1^{\nu_1}x_2^{\nu_2}\prod_{l=1}^M(x_2^q-\la_l
x_1^p)^{n_l},
\end{equation}
with $M\ge 1,$ distinct non-trivial ``roots'' $\la_l\in \bC\setminus\{0\}$ of multiplicities $n_l\in \bN\setminus\{0\},$ and trivial roots of multiplicities $\nu_1,\nu_2\in\bN$ at the coordinate axes.
Here,  $p$ and $q$ have no common divisor, and  $\ka_2/\ka_1=p/q.$

\medskip
More generally,  if $\ka=(\ka_1,\ka_2)$ is any weight with $0<\ka_1\le \ka_2$ such that  the line  $L_\ka:=\{(t_1,t_2)\in\bR^2:\ka_1 t_1+\ka_2t_2=1\}$ is a supporting  line to the Newton polyhedron $\N(\phi) $ of $\phi,$ then  the {\it $\ka$-principal part} of $\phi$
$$
\phi_\ka(x_1,x_2):=\sum_{(\al_1,\al_2)\in L_\ka} c_{\al_1,\al_2} x_1^{\al_1} x_2^{\al_2}
$$
 is a non-trivial polynomial  which is $\ka$-homogeneous of degree $1$ with respect to the dilations associated to this weight as before, and which can be factorized in a similar way as in \eqref{2.4}.  By definition, we then have
\begin{equation}\label{homapprox}
\phi(x_1,x_2)=\phi_\ka(x_1,x_2) +\ \mbox{terms of higher $\ka$-degree}.
\end{equation}

Adaptedness of a given coordinate system can be verified by means of the following criterion (see  \cite{IM-ada}): Denote  by

$$m(\phi_\pr):=\ord_{S^1} \phi_\pr$$
the maximal order of vanishing of $\phi_\pr$ along the unit circle $S^1$ centered at the
origin. The {\it homogeneous distance} of a $\ka$-homogeneous polynomial $P$ (such as $P=\phi_\pr$) is given by $
d_h(P):= 1/{(\ka_1+\ka_2)}=1/|\ka|.$ Notice that $(d_h(P),d_h(P))$ is just the point of intersection of the line given by $\ka_1t_1+\ka_2t_2=1$ with the bi-sectrix $t_1=t_2.$
The height of $P$ can then be computed by means of the formula
\begin{equation}\label{heightp}
h(P)=\max\{m(P), d_h(P)\}.
\end{equation}

 In  \cite{IM-ada} (Corollary 4.3 and  Corollary 2.3), we proved the following characterization of adaptedness of a given coordinate system:
 \begin{proposition}\label{adapted}
The coordinates $x$ are adapted to $\phi$ if and only if one of the following conditions is satisfied:
\medskip

\bee
\item[(a)]  The principal face  $\pi(\phi)$ of the Newton polyhedron  is a compact edge, and $m(\phi_\pr)\le d(\phi).$
\item[(b)] $\pi(\phi)$  is a vertex.
\item[(c)] $\pi(\phi)$ is an unbounded edge.
\ee
\end{proposition}

These conditions had already been introduced by  Varchenko, who  has shown that they are  sufficient for adaptedness when  $\phi$ is analytic. 
 
We also note  that in case (a) we have $h(\phi)=h(\phi_\pr)=d_h(\phi_\pr).$  Moreover, it can be shown that (a) applies whenever $\pi(\phi)$ is a compact edge  and $\ka_2/\ka_1\notin\NN;$   in this case we even have  $m(\phi_\pr)< d(\phi)$ (cf. \cite{IM-ada}, Corollary 2.3).

\medskip
 
 \subsection{Construction of adapted coordinates}\label{consada}

In the case where  the coordinates $\x$ are not adapted to $\phi,$ the previous results show 
  that the principal face $\pi(\phi)$ must be  a  compact edge, that  $m:=\ka_2/\ka_1\in\NN,$   and that  $m(\phi_\pr)> d(\phi).$  One easily verifies that this implies that  $p=m, q=1$ in  \eqref{2.4}, and  that there is at least one, non-trivial  real  root 
  $x_2=\la_lx_1$ of $\phi_\pr$  of multiplicity $n_l=m(\phi_\pr)>d(\phi).$ Indeed, one can show that this root is unique.  Putting $b_1:=\la_l,$ we shall denote the corresponding root  $x_2=b_1 x_1$ of $\phi_\pr$ as its  {\it principal root}.

Changing coordinates 
$$y_1:= x_1, \ y_2:= x_2-b_1x_1^{m},$$
we arrive at a ``better'' coordinate system $y=(y_1,y_2).$  Indeed, this change of coordinates will transform $\phi_\pr$ into a function $\widetilde{\phi_\pr},$  where the principal face of $\widetilde{\phi_\pr}$ will be a horizontal half-line at level $t_2=m(\phi_\pr),$ so that $d(\widetilde{\phi_\pr})>d(\phi),$ and correspondingly one finds that $d(\tilde \phi)>d(\phi),$ if $\tilde\phi$ expresses $\phi$ is the coordinates $y$ (cf. \cite{IM-ada}).  

\medskip
 Somewhat oversimplifying, by iterating this procedure, we essentially arrive at  Var\-chen\-ko's algorithm  for the  construction of  an adapted coordinate system (cf. \cite{IM-ada} for details).

In conclusion, one can show  that there exists a smooth real-valued function $\psi$  (which we may choose as the so-called {\it principal root jet} of $\phi$) of the form
\begin{equation}\label{prjet}
\psi(x_1)=cx_1^{m}+O(x_1^{m+1})
\end{equation}
 with $c\ne 0,$ defined on a neighborhood of the origin  such that an adapted  coordinate system $(y_1,y_2)$ for $\phi$ is given locally near the origin by means of the (in general non-linear) shear
\begin{equation}\label{adaptco}
y_1:= x_1, \ y_2:= x_2-\psi(x_1).
\end{equation}

  In these adapted coordinates, $\phi$ is given by
\begin{equation}\label{phia1}
 \pad(y):=\phi(y_1,y_2+\psi(y_1)).
\end{equation}

\begin{example}\label{ex0}
{\rm 
$$\phi\x:=(x_2-x_1^m)^n+x_1^\ell.$$
Assume that   $\ell>mn.$ Then the coordinates are not adapted. Indeed,  $\phi_\pr\x=(x_2-x_1^m)^n,$ $d(\phi)=1/(1/n+1/(mn))=mn/(m+1)$ and $m(\phi_\pr)=n>d(\phi).$  Adapted coordinates are given by  $y_1:=x_1, y_2:=x_2-x_1^m,$ in which $\phi$ is expressed  by $\pad(y)=y_2^n+y_1^\ell.$
}
\end{example}

  %%%%%%%%%%%%%%%%%%%%%%
%%%%%%%%%%%%%%%%%%%%%%
%%%%%%% Figure 2%%%%%%%%%
%%%%%%%%%%%%%%%%%%%%%%
%%%%%%%%%%%%%%%%%%%%%%
\begin{figure}
\centering
\scalebox{0.35}{\input{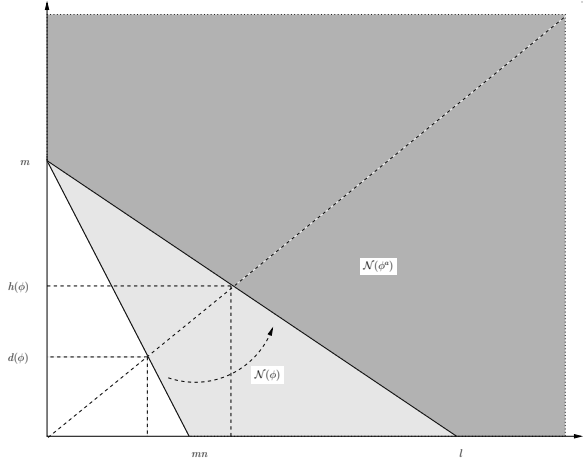}}
\caption{$\phi\x:=(x_2-x_1^{m})^n+x_1^{\ell} \quad (\ell>mn)$}{\label{fig2}}
\end{figure}

 \begin{remark}\label{higherd}
  An alternative proof of  Varchenko's theorem on the existence of adapted coordinates for analytic functions  $\phi$ of two variables has been given by Phong, J. Sturm  and Stein in \cite{PSS}, by means of Pusieux series expansions of the roots of $\phi.$

% (b)In higher dimensions, Varchenko \cite{Va} had already observed that adapted coordinates may not exist. Example ???????
\end{remark}

Let us finally observe that when $m=\ka_2/\ka_1=1$ in the first step of Varchenko's algorithm, then a linear change of coordinates of the form $y_1=x_1, y_2=x_2-b_1x_1$ will transform $\phi$ into a function $\tilde\phi.$  Since all of our problems A - C are invariant under such linear changes of coordinates,  by replacing our original coordinates $\x$ by $\y$ and $\phi$ by $\tilde \phi,$ we may in the sequel always assume without loss of generality that either our coordinates $\x$  are adapted,  or they are not adapted and 
\begin{equation}\label{m1}
m=\ka_2/\ka_1\qquad\mbox{is an integer} \quad \ge 2.
\end{equation}
A linear, non-adapted  coordinate system for which \eqref{m1} holds true will be called {\it linearly adapted} to $\phi.$

 \setcounter{equation}{0}
\section{Problem A: Decay of the Fourier transform of the surface carried measure $d\mu.$ }\label{FT}
Observe first that in view of \eqref{graph1}, we may write $\widehat{\mu}(\xi)$ as an oscillatory integral
$$
\widehat{\mu}(\xi)=:J(\xi)=\int_\Om e^{-i(\xi_3 \phi\x+\xi_1 x_1+\xi_2 x_2)} \eta(x)\, dx,\quad \xi\in\RR^3,
$$
where $\eta\in C_0^\infty(\Om).$ Since $\nabla \phi(0,0)=0,$ the complete phase in this oscillatory integral will have no critical point on the support of $\eta$ unless  $|\xi_1|+|\xi_2|\ll |\xi_3|,$ provided $\Om$ is chosen sufficiently small. Integrations by parts then show that $\widehat{\mu}(\xi)=O(|\xi|^{-N})$ as $|\xi|\to\infty,$ for every $N\in\NN,$ unless $|\xi_1|+|\xi_2|\ll |\xi_3|.$ 

We may thus focus on the latter case. In this case, by writing $\la=-\xi_3$ and $\xi_j=s_j\la, \ j=1,2,$ we are reduced to estimating two-dimensional  oscillatory integrals of the form
$$
I(\la; s):= \int e^{i\la(\phi\x+s_1 x_1+s_2 x_2)} \eta \x \, dx_1\, dx_2,
$$
where we may assume without loss of generality that $\la\gg 1,$ and that $s=(s_1,s_2)\in\RR^2$ is sufficiently small, provided that $\eta$ is supported in a sufficiently small neighborhood of the origin. The complete phase function is thus a small, linear perturbation of the function $\phi.$

If $s=0,$ then the  function $I(\la;0)$ is given by an oscillatory integral of the form
$ \int e^{i\la\phi(x)}\eta(x)\, dx,$  and it is well-known (\cite{bernstein-gelfand}, \cite{atiyah}) that for  any analytic phase function $\phi$  defined on a neighborhood of the origin in $\RR^n$ such that $\phi(0)=0,$ such an  integral  admits an asymptotic expansion as $\la\to\infty$ of the form
\begin{equation}\label{asymp}
\sum_{k=0}^\infty \sum_{j=0}^{n-1} a_{j,k}(\phi)\la^{-r_k}Ê\log(\la)^j,
\end{equation}
provided  the support of $\eta$ is sufficiently small.
 Here, the $r_k$ form an increasing sequence of rational numbers consisting of  a finite number of arithmetic progressions, which  depends only on the zero set of $\phi,$  and the $a_{j,k}$ are distributions with respect to the cut-off function $\eta.$ The proof is based on Hironaka's theorem.  
 
Let us come  back to the case  $n=2.$   Following  \cite{Va} (with as slight modification), we next  define what we like to call {\it Varchenko's   exponent}  $\nu(\phi)\in\{0,1\}:$ 
\smallskip

 If there exists an adapted local coordinate system $y$ near the origin such that the principal face $\pi(\pad)$ of $\phi,$ when expressed by the function $\pad$ in the new coordinates,  is a vertex, and if $h(\phi) \ge 2,$ then we put $\nu(\phi):=1;$ otherwise, we put $\nu(\phi):=0.$ 
 We remark \cite{IM-uniform} that  the first condition is equivalent to the following one:
 {\it If $y$ is any adapted local coordinate system at the origin, then either $\pi(\pad)$ is a vertex, or a compact edge and $m(\pad_\pr)=d(\pad).$}
 
Varchenko \cite{Va}  has shown that the leading exponent in \eqref{asymp} is given by 
$r_0=1/h(\phi),$ and $\nu(\phi)$ is the maximal $j$ for which $a_{j,k}(\phi)\ne 0.$ This implies in particular that 
\begin{equation}\label{estu1}
|I(\la;0)|\le C \la^{-\frac 1{h(\phi)}}\,\log(\la)^{\nu(\phi)}, \quad \la\gg 1 ,
\end{equation}
and this estimate is sharp in  the exponents. Subsequently, Karpushkin \cite{karpushkin} proved  that this estimate is stable under sufficiently small analytic perturbations of $\phi$ (analogous results are known to be wrong in higher dimensions \cite{Va}).  In particular, we find that $J(\la;s)$ satisfies the same estimate \eqref{estu1}  for $|s|$ sufficiently small, so that we obtain the following {\it uniform estimate}  for  $\hat\mu,$ 
\begin{equation}\label{estu2}
|\widehat{\mu}(\xi)|\le C (1+|\xi|)^{-\frac 1{h(\phi)}}\,\log(2+|\xi|)^{\nu(\phi)}, \quad \xi\in\RR^3,
\end{equation}
provided the support of $\rho$ is sufficiently small.

In \cite{IM-uniform}, we proved, by a quite different method,  that Karpushkin's result remains valid for smooth, finite type functions $\phi$, at least for linear perturbations, which led to the following

 \begin{thm}\label{estu3}
Let $S=\mathrm{graph}(\phi)$ be as before, and assume that $\phi$ is smooth and of finite type. Then estimate \eqref{estu2} holds true provided the support of $\rho$ is sufficiently small.
\end{thm}
The special case where $\xi=(0,0,\xi_3)$ is normal to $S$ at the origin  is due to Greenblatt \cite{greenblatt4}.

One can also show   that this estimate is sharp in the exponents even  when $\phi$ is not analytic,  except for the case where the principal face $\pi(\pad)$ is an unbounded edge. Indeed, if  $\pi(\pad)$ is compact, then  \cite{IM-uniform}
$$
I(\la;0) \asymp C \la^{-\frac 1{h(\phi)}}\,\log(\la)^{\nu(\phi)}
$$
as $\la\to+\infty,$ where $C$ is a non-zero constant. However,  when  $\pi(\pad)$ is unbounded then  the following examples, due to  Iosevich and Sawyer \cite{iosevich-sawyer2}, show that there may be a different behavior in general:
Let
 $
\phi(x_1,x_2):=x_2^2+e^{-1/|x_1|^\alpha},
$
with $\al>0;$
then 
$$
|I(\la;0) |\asymp \frac1{\lambda^{1/2}\log{\lambda}^{1/\alpha}}\  \mbox{as}\  \la\to+\infty,
$$
whereas $\nu(\phi)=0.$ 

\medskip

It appears likely that even the full analog of Karpushkin's theorem, i.e., stability of estimate  \eqref{estu1} under arbitrary small, smooth perturbations, holds true in two dimensions.
 \medskip
 \subsection[proofA]{Outline of some main ideas of  the proof}
 
Our proof  of Theorem \ref{estu3} is based on a certain decomposition scheme  for the given surface $S,$ related to the Newton polyhedra of $\phi$ respectively $\pad,$  which has been inspired by the work of Phong and Stein, in combination with re-scaling arguments.  Since the same type of decompositions plays an important role also in the more involved proofs of our results on the other  problems B and C  (which require also further, more  subtle  refinements of them), I shall outline some of the major ideas used in the proof subsequently.

\medskip

We first notice that the case where  $h(\phi)<2$  is covered by  Duistermaat's work \cite{duistermaat}  (notice that Duistermaat proves estimates of the form \eqref{estu2}  without the presence of a logarithmic factor $\log(2+|\xi|),$ even for a wider  class of  phase functions). Note that  the required estimates can also be derived easily from the normal forms in Theorem \ref{normform}.
\medskip

I shall therefore subsequently assume that 
$$
h:=h(\phi)\ge 2.
$$

In many situations, one can reduce the problem to a one-dimensional one by means of van der Corput's  lemma, respectively the following (not quite straight-forward) consequence  of it,  whose formulation goes back to  J.\,E. Bj\"ork  (see \cite{domar}) and  G.\,I.~Arhipov \cite{arhipov}.

\begin{lemma}\label{corput}
Assume that $f$ is a smooth real valued function defined on an interval $I\subset \RR$ which is of polynomial type  $n\ge 2\ (n\in\NN)$, i.e., there are positive constants $c_1,c_2>0$ such that
$$
c_1\le\sum^n_{j=2}|f^{(j)}(s)|\le c_2\quad\mbox{for every}\  s\in I.
$$
Then for $\la\in\RR,$
$$
\Big|\int_{I}e^{i\la f(s)} g(s)\, ds\Big| \le C\| g\|_{C^1(I)} (1+|\la|)^{-1/n},
$$
where the constant $C$ depends only on the constants $c_1$ and $c_2.$
\end{lemma}

\subsubsection{The case where the coordinates are  adapted to $\phi$  }\label{adaptedc}
\medskip

Let us write $d=d(\phi),$  and recall that here $h=d,$  since  the coordinates are adapted.  Assume for instance that the principal face $\pi(\phi)$ is a compact edge. By decomposing $\RR^2$ into the half-spaces  $\RR^2_\pm:=\RR\times \RR_{\pm}$, we may also assume that the integration in $J(\xi)$ takes place over one of these half-spaces only, say, $\RR^2_+.$  
Let $\ka$  be the principal weight, with associated dilations
$\de_r\x=(r^{\ka_1}x_1,r^{\ka_2} x_2).$ Recall also that then $\phi_\pr=\phi_\ka$ is $\de_r$-homogeneous of degree 1.

We fix a suitable smooth cut-off function $\chi$ on $\RR^2$ supported in an annulus $\A$  on which $|x|\sim 1,$ such that the functions $\chi_k:=\chi\circ \de_{2^k}$ form a partition of unity,  and then decompose 
 $$J(\xi)=\sum_{k=k_0}^\infty J_k(\xi),$$
  where
$$
J_k(\xi):=\int_{\RR^2_+}e^{-i(\xi_3\phi(x)+\xi_1x_1+\xi_2x_2)}\eta(x)\chi_k(x)\, dx.
$$
Scaling by $\de_{2^{-k}},$ we see that 
  \begin{equation}\label{intk1}
J_k(\xi)=2^{-k|\ka|} \int_{\RR^2_+}e^{-i\Big(2^{-k}\xi_3\phi^k(x)+2^{-k\ka_1}\xi_1x_1+2^{-k\ka_2}\xi_2x_2\Big)}\eta(\de_{2^{-k}}(x))\chi(x)\, dx,
\end{equation}
with $\phi^k(x):=2^k\phi(\de_{2^{-k}}x).$ Notice that in view of \eqref{homapprox}, 
$$\phi^k(x)=\phi_\ka (x)+ \mbox{error term}.$$

\smallskip

 {\it We claim that given any point $x^0\in \A,$ we can
find a unit vector $e\in\RR^2$ and some $j\in
 \NN$  with $2\le j\le h$ such that $\pa_e^{j}\phi_\ka(x^0)\ne 0,$ where $\pa_e$ denotes the  partial derivative in direction 
 of $e.$ }

Indeed, if $\nabla \phi_\ka(x^0)\neq0,$ then the homogeneity of $\phi_\ka$ and  Euler's homogeneity relation  imply  that $\rank(D^2\phi_\ka(x^0))\ge 1.$   Therefore, we can find a unit vector $e\in\bR^2$ such that $\pa_e^2\phi_{\ka}(x^0)\ne
0.$ 
And, if $\nabla \phi_\ka(x^0)=0,$ then by  Euler's homogeneity relation we have $\phi_\ka(x^0)=0$ as
well. Thus the function $\phi_\ka$ vanishes in $x^0$  of order $j\ge 2.$ But,  in view of Proposition \ref{adapted},  then $j\le m(\phi_\pr)\le d=h,$   which verifies the claim.

\medskip
For $k\ge k_0$ sufficiently large we can thus apply  van der Corput's lemma  to the integration along  lines parallel to the direction $e$ in the integral defining $J_k(\xi)$ near the point $x^0.$ Applying Fubini's theorem and  a partition of unity argument, we thus obtain
\begin{eqnarray}\label{estjk}\nonumber
 |J_k(\xi)|&\le& C\|\eta\|_{C^3(\RR^2)}\, 2^{-k|\ka|}(1+2^{-k}|\xi_3|)^{-1/j}\\
 &\le& C\|\eta\|_{C^3(\RR^2)}\, 2^{-k|\ka|}(1+2^{-k}|\xi|)^{-1/M},
\end{eqnarray}
where $M$ denotes the maximal  $j$ that arises in this context.

Summation in $k$ then yields the following estimates:
\begin{equation}\label{2.3}
|J(\xi)|\le C \|\eta\|_{C^3(\RR^2)}\,\left\{  \begin{array}{cc}
 (1+|\xi|)^{-1/M}, & \mbox{   if } M|\ka|>1\ ,\hfill\\
(1+|\xi|)^{-1/M}\, \log(2+|\xi|), & \mbox{   if } M|\ka|=1\ ,\hfill\\
(1+|\xi|)^{-|\ka|}, & \mbox{   if } M|\ka|<1\
.\hfill
\end{array}\right.
\end{equation}

 However, since we are assuming that $\pi(\phi)$ is a compact edge, we have that  $1/|\ka|=d(\phi)=h,$ and moreover $M\le h.$ This implies $|\ka |M\le 1. $ Since we have seen that here  $\nu(\phi)=1$ if and only if $M=m(\phi_\pr)=h,$ i.e., if and only if $M|\ka|=1,$ we obtain  estimate \eqref{estu2}.

\subsubsection{The case where the coordinates are not adapted to $\phi$  }\label{nonadaptedc}
\medskip

Here, in a {\it first step}, we may reduce to a narrow neighborhood of the principal root. 

Indeed,  away from the principal root of $\phi_\pr,$ we can argue in the same way as before, since the multiplicity of any real root of $\phi_\pr$ different from the principal root is bounded by $d\le h.$ I.e., we can reduce to a narrow $\ka$-homogeneous neighborhood of the curve $x_2=b_1x_1^{m},$ of the form
\begin{equation}\label{restdomain}
|x_2-b_1x_1^{m}|\le \ve x_1^{m},
\end{equation}
say by means of a function $\rho_1(x):=\chi_0((x_2-b_1x_1^m)/(\ve x_1^m)),$ where $\chi_0$ is a suitable smooth bump function supported in the interval $[-1,1]$ and $\ve>0$ is sufficiently small. I.e., in place of $J(\xi),$  it suffices to estimate $J^{\rho_1}(\xi),$ where we write 
$$J^\chi(\xi):=\int_{\RR^2_+}e^{-i(\xi_3\phi\x+\xi_1x_1+\xi_2x_2)}\eta(x)\, \chi(x)\, dx
$$
if $\chi$ is any integrable function.

\medskip
{\it Second Step: Domain decomposition into ``homogeneous'' domains $D_l$ and transition domains $E_l.$}
In oder to study the contribution $J^{\rho_1}(\xi)$  by the domain \eqref{restdomain}, we change to the adapted  coordinates $y$  and essentially  re-write
\begin{equation}\label{form1}
J^{\rho_1}(\xi)=\int_{\RR^2_+}e^{-i(\xi_3\pad\y+\xi_1y_1+\xi_2\psi(y_1)+\xi_2y_2)}\tilde\eta(y)\, \tilde\chi_0\Big(\frac{y_2}{\ve y_1^m}\Big)\, dy,
\end{equation}
where $\tilde\eta$ and $ \tilde\chi_0$ have properties similar to $\eta,$ respectively $\chi_0.$

 Let us  denote the vertices of the Newton polyhedron $\N(\pad)$ by 
 $(A_l,B_l), \  l=0,\dots,n,$ where  we assume that they are ordered so that $A_{l-1}< A_{l},\ l=1,\dots,n,$ with associated compact edges given by the intervals
  $\ga_l:=[(A_{l-1},B_{l-1}), (A_l,B_l)], l=1,\dots,n.$ The unbounded horizontal edge with left endpoint $(A_n,B_n)$ will be denoted by  $\ga_{n+1}.$ To each of these edges $\ga_l,$ we associate the weight
$\ka^l=(\ka^l_1,\ka^l_2),$ so that $\ga_l$ is  contained in  the line
$$L_l:=\{(t_1,t_2)\in \bR^2:\ka^l_1t_1+\ka^l_2 t_2=1\}.$$ 
For $l=n+1,$ we have $\ka_1^{n+1}:=0,\ka_2^{n+1}=1/B_{n}.$ We denote by 
$$a_l:=\frac {\ka^l_2}{\ka^l_1}, \quad l=1,\dots, n$$
the  reciprocal of the slope of the line $L_l.$ For $l=n+1,$ we formally set $a_{n+1}:=\infty$.

If $l\le n,$  the $\ka^l$-principal part $\pad_{\ka^l }$ of $\pad$ corresponding to the supporting line $L_l$  is of the form
\begin{equation}\label{normalhom}
\pad_{\ka^l }(y)=c_l\, y_1^{A_{l-1}} y_2^{B_l}\prod_\al \Big(y_2-c^\al_l y_1^{a_{l}}\Big)^{N_\al}
\end{equation}
(cf. \cite{IKM-max}).

%%%%%%%%%%%%%%%%%%%%%%
%%%%%%%%%%%%%%%%%%%%%%
%%%%%%% Figure 3 %%%%%%%%%
%%%%%%%%%%%%%%%%%%%%%%
%%%%%%%%%%%%%%%%%%%%%%
\begin{figure}
\centering
\vskip-5cm
\scalebox{0.35}{\input{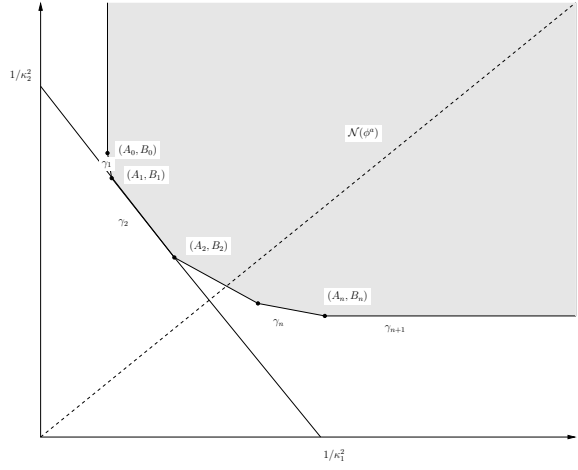}}
\caption{Edges and weights}{\label{fig3}}
\end{figure}

\begin{remark}\label{link}
When $\phi$ is analytic, then this expression is linked to the Puiseux series expansion of roots of $\phi$ as follows \cite{phong-stein} (compare also \cite{IM-ada}):

We may then factorize
$$\pad(y_1,y_2)=U(y_1,y_2)y_1^{\nu_1} y_2^{\nu_2} \prod_{r}(y_2-r(y_1)),
$$
 where the product is indexed by  all non-trivial roots   $r=r(y_1)$ of $\pad$ (which may also be empty)  and where $U(0,0)\ne 0.$ Moreover,  these roots can be expressed  in a small neighborhood of $0$ as   Puiseux series 
$$
r(y_1)=c_{l_1}^{\al_1}y_1^{a_{l_1}}+c_{l_1l_2}^{\al_1\al_2}y_1^{a_{l_1l_2}^{\al_1}}+\cdots+
c_{l_1\cdots l_p}^{\al_1\cdots \al_p}y_1^{a_{l_1\cdots
l_p}^{\al_1\cdots \al_{p-1}}}+\cdots,
$$
where
$$
c_{l_1\cdots l_p}^{\al_1\cdots \al_{p-1}\be}\neq c_{l_1\cdots
l_p}^{\al_1\cdots \al_{p-1}\ga} \quad \mbox{for}\quad \be\neq \ga,
$$

$$
a_{l_1\cdots l_p}^{\al_1\cdots \al_{p-1}}>a_{l_1\cdots
l_{p-1}}^{\al_1\cdots \al_{p-2}},
$$
with strictly positive exponents $a_{l_1\cdots l_p}^{\al_1\cdots \al_{p-1}}>0$ and non-zero complex coefficients $c_{l_1\cdots l_p}^{\al_1\cdots \al_p}\ne 0,$ and where we have kept enough terms to distinguish between all the non-identical roots of $\pad.$ The leading exponents in these series are the numbers
$$
a_1<a_2<\cdots <a_n.
$$

One can therefore group the roots into the clusters of roots $[l], l=1,\dots,n,$ where the $l$'th cluster $[l]$ consistes of all roots with leading exponent $a_l.$

Correspondingly, we can decompose 
$$
\pad(y_1,y_2)=U(y_1,y_2)y_1^{\nu_1} y_2^{\nu_2} \prod_{l=1}^n
\Phi_{[l]}(y_1,y_2),
$$
where
$$
\Phi_{[l]}(y_1,y_2):=\prod_{r\in [l]}(y_2-r(y_1)).
$$
 \end{remark}
 
 Observe the following:   If $\de_s^{l}\x=(s^{\ka_1^{l}}x_1,s^{\ka_2^{l}}x_1), \ s>0,$ denote the dilations associated to the weight $\ka^{l},$ and if $r\in[l_1]$ is a root in the cluster $[l_1],$ then one easily checks that for  $y=\y$ in a bounded set we have $\de_s^{l}y_2=s^{\ka_2^{l}}y_2$ and 
 $r(\de_s^{l}y_1)=s^{a_{l_1}\ka_1^{l}}c_{l_1}^{\al_1}y_1^{a_{l_1}}(1+ O(s^{\ve}))$ as $s\to 0,$ for some $\ve>0.$
  Consequently, 
  \begin{equation}\nonumber
\de_s^{l}y_2-r(\de_s^{l}y_1)=(1+ O(s^{\ve}))\left\{  \begin{array}{cc}
-s^{a_{l_1}\ka_1^{l}}\, c_{l_1}^{\al_1}y_1^{a_{l_1}},& \mbox{   if } l_1<l\ ,\hfill\\
s^{\ka_2^l} \,( y_2-c_l^{\al_l}y_1^{a_l}), & \mbox{   if } l_1=l\ ,\hfill\\
s^{\ka_2^{l}}y_2, & \mbox{   if } l_1>l.\
\hfill
\end{array}\right.
\end{equation}
This shows that the $\ka^l$-principal part of $\pad$ is given by 
 \begin{equation}\label{newprinc}
\pad_{\ka^l}=C_l y_1^{\nu_1+\sum_{l_1<l} |[l_1]| a_{l_1}} y_2^{\nu_2+\sum_{l_1>l} |[l_1]|}\prod_{\al_1}(y_2-c_l^{\al_1}y_1^{a_l})^{N_{l,\al_1}},
\end{equation}
where $N_{l,\al_1}$ denotes the number of roots in the cluster $[l]$ with leading term $c_l^{\al_1}y_1^{a_{l}}.$ A look at the Newton polyhedron reveals that the exponents of $y_1$ and $y_2$ in \eqref{newprinc} can  be expressed in terms of the vertices $(A_j,B_j)$ of the Newton polyhedron: 
$$
\nu_1+\sum_{l_1<l} |[l_1]| a_{l_1}=A_{l-1},\ \nu_2+\sum_{l_1>l} |[l_1]|=B_l.
$$
Notice also that 
$$
\prod_{\al_1}(y_2-c_l^{\al_1}y_1^{a_l})^{N_{l,\al_1}}=\prod_{l=1}^n (\Phi_{[l]})_{\ka^l}.
$$

Comparing this with \eqref{normalhom},  the close relation  between the   Newton polyhedron of $\pad$ and the Pusieux series expansion of roots becomes evident, and accordingly we  say that the edge $\ga_l:=[(A_{l-1},B_{l-1}) ,(A_l,B_l)]$ is
{\it associated to the cluster of roots} $[l].$

\bigskip

Next we  choose the integer $l_0\ge 1$ such that
 $$
a_1<\dots <a_{l_0-1}\le m<a_{l_0}<\dots< a_l<a_{l+1}<\dots<a_n.
$$
Since the original coordinates $x$ were assumed to be non-adapted,  the vertex  $(A_{l_0-1},B_{l_0-1})$ will lie   strictly above the bisectrix, i.e.,  $A_{l_0-1}<B_{l_0-1},$ .

\medskip
Let us consider the case where  the principal face of the Newton polyhedron of $\pad$ is a {\bf compact edge} (the other cases require modified arguments).  We choose $\la>l_0$ so  that the edge $\ga_\la=[(A_{\la-1},B_{\la-1}) ,(A_\la,B_\la)]$ is the principal face $\pi(\pad)$ of the Newton polyhedron of $\pad$ (cf. Figure 3, where $\la=3.$) In the case where $m(\pad_\pr)=d(\pad),$ it can easily be seen that by running Varchenko's algorithm one step further, we can pass to new adapted coordinates in which the principal face is a vertex.

We therefore may assume that $m(\pad_\pr)<d(\pad),$  so that $\nu(\phi)=0.$
\medskip

We shall   narrow down the domain \eqref{restdomain} to a neighborhood  $D_\la$ of  the principal root jet of the form
\begin{equation}\label{restdomain2}
 |x_2-\psi(x_1)|\le N_\la  x_1^{a_\la},
\end{equation}
where $N_\la$ is a  constant to be chosen later. This domain is
 $\ka^\la$-homogeneous in the adapted coordinates $y.$

%%%%%%%%%%%%%%%%%%%%%%
%%%%%%%%%%%%%%%%%%%%%%
%%%%%%% Figure 4 %%%%%%%%%
%%%%%%%%%%%%%%%%%%%%%%
%%%%%%%%%%%%%%%%%%%%%%
%\begin{figure}
%\centering
%\scalebox{0.35}{\input{roots.pstex_t}}
%\caption{Clusters of roots}{\label{fig4}}
%\end{figure}
\medskip
To this end, following an idea from \cite{phong-stein}, we decompose the difference set of the  domains \eqref{restdomain} and \eqref{restdomain2}  into   the domains
 $$
D_l:=\{\x:\ve_l  x_1^{a_l}< |x_2-\psi(x_1)|\le N_l x_1^{a_l}\},\quad l=l_0,\dots,\la-1,
$$
 and  the intermediate  domains
$$E_l:=\{\x:N_{l+1} x_1^{a_{l+1}}<|x_2-\psi(x_1)| \le \ve_l x_1^{a_l}\}, \quad l=l_0,\dots,\la-1,
$$
as well as 
$E_{l_0-1}:=\{\x:N_{l_0} x_1^{a_{l_0}}<|x_2-\psi(x_1)| \le \ve_1 x_1^{m}\}.
$
Here,  the $\ve_l>0 $ are  small and the $N_l>0$ are  large parameters to be chosen suitably.

Observe that the  domain 
$$
D^ a_l:=\{\y:\ve_l  y_1^{a_l}< |y_2|\le N_l y_1^{a_l}\}
$$
corresponding to $D_l$ is $\ka^l$-homogeneous in the adapted coordinates $y$ given by \eqref{2.4},  and contains the cluster of roots $[l],$ if $\ve_l$ and $N_l$ are chosen sufficiently small, respectively  large, while the domain $E_l^a$ corresponding to $E_l$ can be viewed as a domain of transition between two different homogeneities.

\medskip
 The contributions by the domain $D_l$ to $J^{\rho_1}(\xi)$  can  again be estimated in a similar way as we did in the adapted case, by using dyadic decompositions and subsequent re-scalings by means of  the dilations $\de^l_r$ associated to the weight $\ka^l.$ 

More precisely, the corresponding $k$'th term  will be given by the following analogue of \eqref{intk1}
$$
J_k(\xi)=2^{-k|\ka^l|} \int_{\RR^2_+}e^{-i\Big(2^{-k}\xi_3\phi^k(y)+ \xi_2\psi(2^{-k\ka^l_1}y_1)+2^{-k\ka^l_1}\xi_1y_1+2^{-k\ka^l_2}\xi_2y_2\Big)}\eta(\de^l_{2^{-k}}y)\chi(y)\, dy,
$$
where $\phi^k(y)=\pad_{\ka^l} (y)+ \mbox{error term}.$

Notice, however, that since $1-m\ka_1^l>\ka^l_2-m\ka_1^l>0,$  the contribution of the non-linearity $\psi$ to the complete phase of the corresponding oscillatory integrals may be   large, compared to the term containing $\phi^k,$  so that we are only allowed to apply van der Corput's estimate  to the integration with respect to the variable $y_2$   if we want to reduce to one-dimensional oscillatory integrals!
\medskip

  We therefore need a {\it control on the multiplicities of roots of $\pa_2^2\pad_{\ka^l}$} at points $y^0$ in the corresponding annulus $\A$ not lying on the $y_1$ axis (note that the latter  corresponds to the principal root jet in the coordinates $y$).  Indeed, it can be shown (cf.  Propostion 2.3 (b) in \cite{IKM-max}) that these multiplicities are bounded by $d_h(\pad _{\ka^l})-2,$ where $d_h(\pad _{\ka^l})$ denotes the homogeneous distance of $\pad _{\ka^l},$ and it is evident from the geometry of the Newton polyhedron of $\pad$ that 
$d_h(\pad _{\ka^l})<d(\pad)=h,$ so that for every point $y^0$ in $\A\cap D_l$ there is some $j\in\{2,\dots, h-1\}$ such that  
$$\pa_2^j\pad_{\ka^l}(y^0)\ne 0.
$$

\medskip
As for the contributions by the domains $E_l,$ here we perform a separate dyadic decomposition in both variables $y_1$ and $y_2,$ so that we geometrically decompose $E_l$ into dyadic rectangles of size $2^{-j}\times 2^{-k},$ and then re-scale in both variables so that these rectangles become the standard cube, say, $[1,2]\times [1,2].$

Now, the estimates of   Section 11 in \cite{IKM-max} show that the  phase functions $\pad_{j,k}$ that one obtains after these re-scalings  satisfy the estimate
$$\pa_2^2\pad_{j,k}(y^0)\ne 0 \ \mbox{ for every} \ y^0\in [1,2]\times [1,2].$$ 
Since $h\ge 2,$ this clearly suffices to obtain the necessary order of decay of the Fourier transform of these dyadic pieces. Moreover, scaling back to the original dyadic rectangles, a careful analysis of the dependency of the corresponding estimates on the parameters $j,k$  shows that it is indeed possible to sum theses estimates and obtain the same type of estimate for the contributions by the domains  $E_l$ as for the  domains $D_l,$ even without logarithmic factor. 
\smallskip

{\it Note that so far, we have always been able to reduce our estimations to van der Corput's lemma, respectively, Lemma \ref{corput}, i.e., to one-dimensional oscillatory integrals.}

\subsubsection{Third Step: Study of the  contribution by the homogenous domain $D_\la$ containing the  principal root jet}\label{near}

What remains to be estimated is the contribution of the  domain  \eqref{restdomain2} to $J(\xi).$ 
As we shall see, the study of this domain will in certain cases require the estimation of genuinely 2-dimensional oscillatory integrals, and a reduction to the one-dimensional case is no longer possible.

We recall that according to our convention
\begin{equation}\label{convb}
m(\pad_\pr)<d(\pad)=h,
\end{equation}
so that $\nu(\phi)=0.$ 
In the adapted coordinates $y,$ the domain  \eqref{restdomain2} is given by 
$
 |y_2|\le N_\la  y_1^{a_\la},
$
and we can cover it by a finite number of  $\ka^\la$-homogeneous subdomains of the form
$
 |y_2-cy_1^{a_\la}|\le \ve_0  y_1^{a_\la},
$
where $c\in[-N_\la,N_\la],$ and where, for a given $c,$ we may choose $\ve_0>0$ suitably small. 

Writing $\psi(x_1)=x_1^{m}\om(x_1),$ with a smooth function $\om$ satisfying  $\om(0)\ne 0,$  we can thus reduce to estimating oscillatory integrals of the form
\begin{equation}\label{7.3}
J^c(\xi)=\int_{\RR_+^2} e^{iF(y,\xi)}\rho\Big(\frac{ y_2-cy_1^{a_\la}}{\ve_0 y_1^{a_\la}}\Big)\eta(y)\,dy,
\end{equation}
with a phase function
$$
F(y,\xi):=\xi_3\pad(y)+\xi_1y_1+\xi_2 y_1^m\om(y_1)+\xi_2y_2
$$
depending on $\xi\in\RR^3.$

Arguing in a similar way as in the case of adapted coordinates, and recalling that $\pad_\pr= \pad_{\ka^\la},$ we may again perform a dyadic decomposition and re-scale  by means of the dilations $\de^\la_r,$ in order to write

$$
J^c(\xi)=\sum_{k=k_0}^\infty J_k(\xi),
$$
where
\begin{equation}\label{jk}
J_k(\xi)=2^{-|\ka^\la|k}\int e^{i2^{-k}\xi_3F_k(y,s)} \rho\Big(\frac{ y_2-cy_1^{a_\la}}{\ve_0y_1^{a_\la}}\Big)\eta(\de^\la_{2^{-k}}y)\,\chi(y)\,dy,
\end{equation}
where
$$
F_k(y,s):=\pad_\pr\y+s_1y_1+S_2y_1^{m}\om(2^{-\ka^\la_1k}y_1)+s_2y_2+{\rm error},
$$
where $s:=(s_1,s_2,S_2)$  is given by 
$$
s_1:=2^{(1-\ka^\la_1)k}\frac{\xi_1}{\xi_3},\ s_2:=2^{(1-\ka^\la_2)k}\frac{\xi_2}{\xi_3},\
S_2:=2^{(\ka^\la_2-m\ka^\la_1)k}s_2.
$$

%Here, we have decomposed $\pad=\pad_{\pr}+\phi_r=\pad_{\ka^\la}+\phi_r,$ where $\phi_r$ consists of terms of $\ka^\la$-degree strictly bigger than $1$  and can thus be considered as error term.

Note that $2\le m<a_\la=\ka^\la_2/\ka^\la_1$ and $k\gg  1,$ so that $|S_2|\gg  |s_2|,$ and that $$
y_1\sim 1 \ \mbox{and} \ |y_2-cy_1^{a_\la}|\lesssim \ve_0
$$
for $y$ in the support of the integrand of $J_k(\xi).$
Recall also   that we are  assuming  that $ \ |\xi|\sim |\xi_3|.$

One is thus led to the   {\it estimation of   oscillatory integrals depending on  certain parameters} (here $s_1,s_2,S_2$) which may have various relative sizes. 

\medskip
The case where $|S_2|  \ge M$ for some sufficiently  large constant $M\gg 1$ can easily be treated by means of  van der Corput's lemma  applied to the $y_1$- integration, so let us assume that  $|S_2|  < M.$    Then $|s_2|\ll1,$ provided we have chosen $k_0$ sufficiently large.

We  may also easily reduce to the case where 
\begin{equation}\label{2.15}
\partial_2^j\pad_{\pr}(1,c)=\partial_2^j\pad_{\ka^\la}(1,c)= 0 \ \mbox{for}\ 1\le j<h,
\end{equation}
for otherwise an integration by parts in $y_2$ (if $j=1$) or a simple application of Lemma \ref{corput} yields a suitable estimate as before.

The case where $c>0$ can easily be reduced to the case $c=0$ by  performing  another change of variables $y_2\mapsto y_2+cy_1^{a_\la}$ in the integral defining $J_k(\xi).$ Indeed,  one can show that  our assumption  \eqref{2.15} implies that necessarily  $a_\la=\ka^\la_2/\ka^\la_1\in\NN$ (cf. Corollary 3.2 (iii)  in \cite{IKM-max}),  and one checks that the new coordinates are again adapted to $\phi.$ 
 \medskip

 So, let us assume that $c=0.$
Then necessarily  $\pad_\pr(1,0)\ne  0,$ for otherwise  $\pad_\pr$ would have  a root of multiplicity  at least $h$ at $(1,0),$ which would contradict \eqref{convb}.

Assuming  without loss of generality that   $\pad_\pr(1,0)=1,$
we  can then write (compare \cite{IKM-max}, Subsection 9.1)
$$
\pad_\pr\y=y_2^BQ\y+y_1^n,
$$
where $Q$ is a $\ka^\la$-homogeneous polynomial such that
$Q(1,0)\ne 0,$ and where $B\ge h>2.$

\bigskip
Recall that  $|S_2|<M,$  so that  $|s_2|\ll1.$  Moreover, the case where  $|s_1|  \ge N$ for some sufficiently  large constant $N\gg 1$ can easily be dealt with by means of  an integration by parts in $y_1,$ so let us assume that $|s_1| < N.$

It will then suffice to show that, given any point $(s_1^0,S_2^0)\in   [-M,M]\times [-N,N]$ and any point $y_1^0\sim 1,$ there exist a  neighborhood $U$ of
$(s_1^0,S_2^0),$ a neighborhood $V$ of $y_1^0$ and some $\si >1/h$ such that we have an estimate of the form
\begin{equation}\label{2.17}
|J_k(\xi)|\lesssim\frac{2^{-k|\ka|}}{(1+2^{-k}|\xi|)^{\si}}
\end{equation}
for every $(s_1,S_2)\in U,$ provided the function $\chi$ in the definition of $J_k(\xi)$ is supported in $V$ and $\ve_0$ and $k$ are chosen sufficiently small, respectively large.     Summing  over all $k,$ this will clearly imply an estimate as in \eqref{estu2}, even without logarithmic factor.

\medskip
To this end, first notice that for $(s_1,S_2)\in U$  and $k$ sufficiently large, the function  $F_k(y,s)$ can be viewed as a small $C^\infty$- perturbation of the function
$$
F_{\pr}(y):=y_2^BQ\y+s^0_1y_1+S^0_2 \om(0)y_1^m+y_1^n.
$$
Thus, if $\nabla F_\pr(y_1^0,0)\ne 0,$ then we obtain \eqref{2.17}, with $\si=1,$ simply by integration by parts.

\medskip
Assume therefore that $(y_1^0,0)$ is a critical point of $F_\pr.$ Then $y_1^0$ is a critical point of  the polynomial function
$$
g(y_1):= s_1^0y_1+S_2^0 \om(0)y_1^m+y_1^n,
$$
 which comprises all terms of $F_\pr$ depending on the variable $y_1$ only. Note that the exponents satisfy $1< m<n,$ since
 $n=1/{\ka^\la_1}>\ka^\la_2/{\ka^\la_1}>m.$  It is then easy to see that $g''$ and $g'''$ cannot vanish simultaneously at the given point $y_1^0,$ so that
  there are positive constants $c_1,c_2>0$ and a compact neighborhood $V$ of $(y_1^0,0)$ such that
$$
c_1\le\sum^3_{j=2}|g^{(j)}(y_1)|\le c_2\quad\mbox{for every}\  y_1\in V.
$$
This implies an analogous estimate for the partial derivatives $\partial_{1}^jF_k(y_1,y_2,s)$ of $F_k,$ uniformly for $(s_1,S_2)\in U$   and the range of $y_2$ that we consider, provided we choose $U$ and $\ve_0$ sufficiently small.
Applying the van der Corput type estimate in Lemma \ref{corput}, we thus obtain   estimate \eqref{2.17}, with $\si=1/3,$  so that we are done provided $h>3.$ Notice also that if $g''(y_1^0)\ne 0,$  then by  the same type of argument we see that \eqref{2.17} holds true with $\si=1/2>1/h.$

\medskip

We may thus  finally  assume that $2<h \le 3,$ and that $g'(y_1^0)=g''(y_1^0)=0.$  In this case we have
$$
\frac1{\ka^\la_1+\ka^\la_2}=h\le3 \quad \mbox {and}\quad \frac{\ka^\la_2}{\ka^\la_1}> m\ge2,
$$
so  that  $1/\ka^\la_2< 9/2.$

Note that $B\le 1/\ka^\la_2$ is a positive integer, and   $h\le B<9/2,$ so that either
 $B=4$ or $B=3$. We translate the critical point $(y_1^0,0)$ of $F_\pr$ to the origin by considering the function
 $$
 F_\pr^\sharp(z):=F_\pr(y_1^0+z_1,z_2)-g(y_1^0)=z_2^B Q(y_1^0+z_1,z_2)+
 \frac 16 g^{(3)}(y_1^0)\, z_1^3+\dots .
 $$
It is easy to see that this function has height $h^\sharp:=h(F_\pr^\sharp)$ given by
$h^\sharp=\frac 1{1/3+1/B}<2.$  
 We can  therefore again apply
 Duistermaat's results in \cite{duistermaat} to  the  corresponding part of the oscillatory integral $J_k(\xi)$ and obtain  estimate \eqref{2.17}, with $\si=1/h^\sharp>1/h.$ Note here that the estimates in \cite{duistermaat} are stable under small perturbations.

\setcounter{equation}{0}
\section{Problem B: $L^p$- boundedness of the maximal operator $\M$ associated to the hypersurface  $S.$}\label{MAX}

We next turn to Problem B. As has already been mentioned in  the introduction, the first, fundamental result on this problem is due to Stein \cite{stein-sphere}, who proved that for $n\ge 3,$ the spherical maximal function is bounded  on $L^p(\RR^n)$  for every $p> n/(n-1).$ The analogous estimate in dimension  $n=2$ turned to require even deeper methods and was later established  by Bourgain \cite{bourgain-circle}. These results became  the starting point for intensive  studies of various classes of maximal operators associated to subvarieties. Stein's  monograph \cite{stein-book} is an excellent reference to many of these developments. From these early works,  the influence of geometric properties of $S$  on the validity of $L^p$-estimates for the maximal operator $\M$ became evident.
For instance, Greenleaf \cite{greenleaf} proved that $\M$ is bounded on $L^p(\RR^n)$ if $n\geq 3$ and $p>{n}/{(n-1)},$
provided $S$ has everywhere non-vanishing Gaussian curvature. 
In contrast, the case where the Gaussian curvature vanishes at some point is still wide open, and complete answers are  available at present only for  finite type curves in the plane  \cite{iosevich-curves}, and finite type hypersurfaces in $\RR^3$ and $p>2$ (cf. Theorem \ref{maxs1.2}). 

\medskip
Let me come back to our hypersurface $S$ in $\RR^3.$
Recall  from the introduction that we assume that $S$ satisfies the transversality condition of \S \ref{intro},  so  that, by localizing to a small neighborhood of a given point $x^0\in S$ and applying a  suitable linear change of coordinates, we may assume that $x^0=0$ and that  $S$ is    given  as the graph $S={\rm graph}(1+\phi),$ where $\phi$ is a smooth, finite type function defined on a neighborhood of the origin and satisfying $\phi(0,0)=0, \nabla\phi(0,0)=0.$
We then define the height of $S$ at $x^0$ by 
$$h(x^0,S):=h(\phi).$$
 It easily seen that this  notion is  invariant under affine linear changes of coordinates in the ambient space  $\RR^3.$
Recall  also that $\M f(x):=\sup_{t>0}|A_tf(x)|,$  where $A_t$ denotes the averaging operator over the $t$-dilate of $S$ given by 
$$
A_tf(x):=\int_{S} f(x-ty) \rho(y) \,d\si(y),\quad t>0.
$$
We can then state our  main result from  \cite{IKM-max} (Theorems 1.2, 1.3), which gives an almost complete answer to the question of $L^p$-boundedness of the maximal operator $\M$ when $p>2:$

\begin{thm}\label{maxs1.2} Assume the hypersurface $S$ satisfies the transversality Assumption \ref{s1.1}.
 
  \begin{itemize}
\item[(i)] If the measure $\rho d\si$ is supported in a sufficiently small neighborhood of $x^0,$ then 
 $\M$ is bounded on $L^p(\RR^3)$ whenever $p>\max\{h(x^0,S),2\}.$
\item[(ii)] If $\M$ is bounded on $L^p(\RR^3) $ for some  $p>1,$ 
and if  $\rho(x^0)>0,$ then 
$p\ge h(x^0,S).$ Moreover, if $S$ is analytic at $x^0,$
then $p>h(x^0,S).$
\end{itemize}
 \end{thm}

\medskip
 \subsection[related]{Related quantities: contact index and sublevel growth}
 In \cite{iosevich-sawyer1}, Iosevich and Sawyer had discovered a very interesting connection between the behavior of maximal functions $\M$ associated to hypersurfaces  $S$ (in $\RR^n$) and  an integrability  index associated to the hypersurface. In order to describe this, assume that   $d_H(x):=\dist(H,x)$  denotes the distance from a point $x$ on $S$ to  a given hyperplane    $H.$  In particular,  if $x^0\in S,$ then $d_{T,x^0}(x):=\dist (x^0+T_{x^0}S, x)$ will denote  the distance from $x\in S$ to the affine tangent plane to $S$ at the  point $x^0.$ The following result has been proved in \cite{iosevich-sawyer1} in arbitrary dimensions $n\ge 2$  and without requiring Assumption \ref{s1.1}.

\begin{namedthm}[Iosevich-Sawyer]\label{maxs1.6}
If the maximal operator $\M$ is bounded on $L^p(\RR^n),$ where $p>1,$ then
\begin{equation}\label{max4.1}
\int_S d_H(x)^{-1/p}\,\rho(x)\,  d\si(x)<\infty
\end{equation}
for every affine hyperplane $H$ in $\RR^n$ which does not pass through the origin.
\end{namedthm}

Moreover, it was  conjectured in \cite{iosevich-sawyer1} Êthat for $p>2$ the condition \eqref{max4.1} is indeed necessary  and sufficient for the boundedness of  the maximal operator $\M$ on $L^p,$ at least if for instance $S$ is compact and $\rho>0.$
\begin{remark}\label{maxs1.7}
Notice that  condition \eqref{max4.1} is easily seen to be true for every affine hyperplane $H$ which is nowhere tangential to $S,$ so that it is in fact  a condition on affine tangent hyperplanes to $S$ only. Moreover, if Assumption \ref{s1.1} is satisfied, then there are no affine tangent hyperplanes which pass through the origin, so that in this case it is a condition on all affine tangent hyperplanes.
\end{remark}

Moreover, it is not very hard to prove (cf. \cite{IKM-max}) that if  $S$ is a smooth hypersurface of finite type in $\RR^3,$ then, for every $p<h(x^0,S),$ 
\begin{equation}\label{max4.2}
\int_{S\cap U} d_{T,x^0}(x)^{-1/p}\,d\si(x)=\infty \quad \mbox{ for every } \ p<h(x^0,S)
\end{equation}
and every  neighborhood $U$ of $x^0.$ And, if $S$ is analytic near $x^0,$  then 
\eqref{max4.2}  holds true  also for $p=h(x^0,S).$

Notice that  this result does not require Assumption \ref{s1.1}.
As an immediate consequence of Theorem \ref{maxs1.2}, Theorem \ref{maxs1.6} and  \eqref{max4.2}  we obtain
\begin{cor}\label{maxs1.9}
Assume that $S\subset \RR^3$ is of finite type and satisfies Assumption \ref{s1.1}, and let $x^0\in S$ be a fixed point. Moreover, let $p>2.$ 

Then, if $S$ is analytic near $x^0,$ 
there exists a neighborhood $U\subset S $ of the point $x^0$ such
that for any $\rho\in C_0^\infty(U)$ with $\rho(x^0)>0$ the associated maximal
operator $\M$ is bounded on $L^p(\RR^3)$ if and only if condition \eqref{max4.1} holds for every affine hyperplane $H$ in $\RR^3$ which does not pass through the origin. 

If $S$ is only assumed to be smooth near $x^0,$ then the same conclusion holds true, with the possible exception of the exponent $p=h(x^0,S).$ 
\end{cor}
This confirms the conjecture by Iosevich and Sawyer in our setting for analytic  $S$, and for smooth, finite-type  $S$ with the possible exception of the exponent $p=h(x^0,S).$ For the critical exponent $p=h(x^0,S),$ if $S$ is not analytic near $x^0,$ examples show that unlike in the analytic case it may happen that $\M$ is bounded on $L^{h(x^0,S)}$ (see, e.g., \cite{iosevich-sawyer2}), and  the conjecture remains open for this value of $p.$

In view of these results, it is natural to define the  {\it uniform contact index} $\ga_u(x^0,S)$ of the hypersurface $S$ at the point $x^0\in S$ as the supremum over the set of all $\ga$  for which   there exists an open  neighborhood $U$ of $x^0$ in $S$  such   that the  estimate
$$
\int_{U\cap S} d_H(x)^{-\ga}\,  d\si(x)<\infty
$$
holds true for every affine hyperplane $H$ in $\RR^n.$ If we restrict ourselves in this definition to the affine  tangent hyperplane $H=x^0+T_{x^0}S$ at the point $x^0,$ we shall call the corresponding index 
the {\it  contact index} $\ga(x^0,S)$ of the hypersurface $S$ at the point $x^0\in S.$  Here, we shall always assume that $\rho(x^0)\ne 0.$ Note that if we change coordinates so that $x^0=0$ and $S$ is  the graph of $\phi$  near the origin, where $\phi$ satisfies \eqref{null}, then the contact index is just the supremum over all $\ga$ such that there is some neighborhood $U$ of the origin so that 
$|\phi|^{-\ga}\in  L^1(U),$ so that the contact index agrees with the {\it critical integrability index of $\phi$} as defined for instance in \cite{cgp}.

A closely related quantity is the {\it sublevel growth rate $\si(\phi),$} defined as the supremum over all $\si>0$ such that there is some constant $C_\si>0$ so that 
$$
\Big|\{x\in U:|\phi(x)|<\ve\}\Big|\le C_\si\ve^\si \quad \mbox{for every }\ \ve>0.
$$
Indeed, it can easily be shown by means of Tchebychev's inequality (cf. \cite{cgp}) that
$$
\si(\phi)=\ga(x^0,S), \quad\mbox{if } x^0=0 \ \mbox{and}\  S={\rm graph}(\phi).
$$

In analogy with Arnol'd's notion of  ``singularity index'' \cite{arnold}, let us finally introduce  the {\it uniform oscillation index} $\be_u(x^0,S)$ of the hypersurface $S$ at the point $x^0\in S$ as the supremum over the set of  all $\be$ such that 
$$
|\widehat{\rho d\si}(\xi)|\le C_\be\,(1+|\xi|)^{-\be},
$$
 for all $\rho$ supported in a sufficiently small neighborhood of $x^0.$ 
If we restrict directions $\xi$  to the normal to $S$ at $x^0,$ then the corresponding decay rate will be called the {\it oscillation index} $\be(x^0,S)$ at $x^0.$

Combining our results with results from \cite{PSS} (compare also \cite{greenblatt3}), we easily obtain the following result for smooth, finite type hypersurfaces $S$ in $\RR^3:$
\begin{cor}\label{maxs4.5}
Let $x^0\in S\subset\RR^3$ be a given point so that $\rho(x^0)>0.$  Then 
 \begin{eqnarray*}
\be_u(x^0,S)= \be(x^0,S)
=\ga_u(x^0,S)= \ga(x^0,S)=1/h(x^0,S).
 \end{eqnarray*}
\end{cor}

\begin{remarks}\label{maxs4.6}
\begin{itemize}
\item [(a)] In dimension $n\ge3,$ the corresponding identities may fail to be true.
Indeed, Varchenko had already observed in \cite{Va} (Example 3)  that  for  graphs  $S$ of functions $\phi$ of three variables the oscillation index $\beta(0,S)$ may differ from $1/h(0,S).$ Moreover, it is easy to give examples of functions $\phi$ of three variables where the contact index $\ga=\ga(0,S)$ is strictly smaller than the oscillation index $\beta=\beta(0,S).$ For instance, this applies to $\phi(x_1,x_2,x_3):= x_3^2-(x_1^2+x_2^2),$ where the method of stationary phase shows that $\be=3/2,$ whereas the factorization $\phi(x_1,x_2,x_3)= (x_3-\sqrt{x_1^2+x_2^2})(x_3+\sqrt{x_1^2+x_2^2})$ shows that $\ga=1$ (compare counterexample 8.1 in \cite{cgp}).

However, the proof of Theorem 1.6 in \cite{greenblatt1} shows that the oscillation index $\be$ can only be different from the contact index  $\ga$ when $1/\ga$ is an odd integer.

To the best of my knowledge, it is not known whether the identities $\be_u(x^0,S)= \be(x^0,S)$ and $\ga_u(x^0,S)= \ga(x^0,S)$ may persist in higher dimensions. Indeed, this question represents a special case of a question posed  by Arnol'd, namely  whether the oscillation index of  a given function $\phi$ is semicontinuos in the following sense:

Suppose $\phi(x,s)$ is a phase function which depends on $x$ in a small neighborhood of the origin in $\RR^n$ and some  small parameters $s\in\RR^k.$ Is there is a small neighborhood of $(0,0)\in\RR^n\times\RR^k$ so that the oscillation index of $\phi(x,s)$ 
at $x^0$ is greater or equal to that  of $\phi(x,0)$ at the origin for every  $(x^0,s^0)$ in this neighborhood?  That question had been answered to the negative in dimensions $n\ge 3$  by Varchenko \cite{Va}, and to the positive for analytic functions depending on two variables by Karpushkin \cite{karpushkin}.
 
 For linear perturbations of a given function $\phi,$ our Theorem \ref{estu3} shows that Arnol'd's conjecture holds true even for smooth, finite type functions $\phi,$ and it seems open whether this kind of stability may still hold true even in higher dimensions, since all  counter examples  to Arnol'd's conjecture which have been found hitherto have non-linear perturbation terms. 
 
 Further important papers dealing with this stability question for general perturbations are, e.g.,  \cite{ccw}, \cite{PSS}.
 
 \item [(b)] Greenblatt \cite{greenblatt-resol1},\cite{greenblatt-resol2},  and independently  Collins, Greenleaf and Paramanik \cite{cgp},  have devised (quite different) algorithms of resolution of singularities which in principle allow to compute the contact index also in higher dimensions. 

 \item [(c)]   If $p\le 2,$  then examples (see, e.g., \cite{io-sa-seeger}) show that neither the notion of height nor that of contact index will  determine the range of exponents $p$ for which the maximal operator $\M$ is $L^p$-bounded. Our work (in progress)  on this case seems to indicate a certain conjecture how to express  this range in terms of Newton polyhedra. Moreover,  for certain surfaces this conjecture  relates to  fundamental open problems in Fourier analysis, such as the  conjectured reverse square function estimate for the cone multiplier (see, e.g.,  \cite{mockenhaupt}, \cite{bourgain-cone}).  
 \end{itemize}
\end{remarks}

\medskip
 \subsection[proofA]{ A few hints on  the proof of Theorem \ref{maxs1.2}}
 
 Since the proof in \cite{IKM-max} is quite involved, I shall here just try to indicate some of the main steps, grossly oversimplifying some of the arguments.
 
 In the first part of the proof, we basically follow the first two steps of the scheme that has been outlined in the preceding paragraph. Recall that in these steps, it had effectively been possible to reduce the problem to a one-dimensional one. The same applies basically to the estimation of the maximal operator $\M,$  because in the first two steps, we have had a sufficiently  good control on the multiplicity of roots by the height $h.$

 In these arguments, the following result plays an analogous role for problem B as the van der Corput type Lemma \ref{corput} played for problem A.

Let $U$ be an open neighborhood of the point $x^0 \in \RR^2,$ and let
$\phi_\pr\in C^\infty(U,\RR)$ such that 
\begin{equation}\label{max4.3}
\pa_2^m\phi_\pr(x^0_1,x^0_2)\neq 0,
\end{equation}
where $m\ge 2.$ Let 
$$\phi=\phi_\pr+\phi_r,$$
where $\phi_r\in C^\infty(U,\RR)$  is a sufficiently small perturbation.
 Denote by 
$S_\ve$ the surface  in $\bR^3$ given  by
$S_\ve:=\{(x_1,x_2,1+\ve \phi(x_1,x_2)): (x_1,x_2)\in U\},$  with $\ve>0,$
and
 consider the averaging operators
$$
A^\ve_tf(x):=\int_{S_\ve} f(x-ty) \psi(y)\,d\si(y),
$$
where $d\si$ denotes the surface measure and  $\psi\in C^\infty_0(S_\ve)$ is
a non-negative cut-off function.  Define the  associated maximal operator by 
$$
\M^\ve f(x):=\sup_{t>0}|A^\ve_tf(x)|.
$$

\begin{prop}\label{maxs4.7}
 Assume that $\phi_\pr$ satisfies \eqref{max4.3} and that the neighborhood $U$ of the
point $x^0$  is sufficiently small. Then  there
exist numbers $M\in\NN$, $\delta>0$, such that for every $\phi_r\in
C^{\infty}(U,\RR)$ with $\| \phi_r\|_{C^M}<\delta$ and every 
$p>m$  there exists a positive constant
$C_p$ such that for $\varepsilon>0 $ sufficiently small the maximal operator $\M^\ve$
satisfies  the following a priori estimate:
\begin{equation}\label{max4.4}
\|{\cal M}^{\varepsilon}f\|_p\le
C_p\, \varepsilon^{-1/p}\| f\|_p,\quad f\in{\cal S}(\RR^3)\, .
\end{equation}
\end{prop}

  This result can  be reduced to the one-dimensional case, i.e.,   the study of maximal functions associated to curves in the plane. Indeed, just consider the "fan" of all 
hyperplanes passing through the $x_1$-axis. This fan will fibre the  given surface into a family of curves, and our dilations leave each of these planes invariant.  

The proof for the analogous result for plane curves then basically follows Iosevich's approach in \cite{iosevich-curves}. 
The case where $m=2$ is indeed the most difficult on; the cases where $m\ge 3$ can then easily be reduced to this case by means of dyadic decompositions and re-scalings.  Note that if $m=2,$ then the corresponding maximal operator in the plane  behaves in a similar way as  the ``circular'' maximal function studied by Bourgain \cite{bourgain-circle}. An alternative approach to Bourgain's  one had later been given by Mockenhaupt, Seeger and Sogge \cite {mockenhaupt-seeger-sogge},  \cite {mockenhaupt-seeger-sogge2}, based on suitable local smoothing estimates for classes of Fourier integral operators. This approach  is stable under small perturbations of the given curve, which is exactly what is needed for our purposes. 

Indeed, in our applications of Proposition \ref{maxs4.7}, we even need to replace $\phi\x$ by $\phi(x_1,x_2-\psi_\ve(x_1)),$ where the function $\psi_\ve$ may blow up like $O(\ve^{-\de})$ for some $\de\in [0,1[,$ and it turns out that this is still admissible, i.e., the estimate \eqref{max4.4} remains valid also in this case.

\medskip
Now, if the coordinates are {\it adapted} to $\phi$ (a similar argument will apply to the  first  step when the  coordinates are not adapted),  then, in analogy with the decomposition  $J(\xi)=\sum_{k=k_0}^\infty J_k(\xi),$ we can dyadically decompose  the surface $S$ as in Subsection \ref{adaptedc}, and accordingly decompose 
  $$
A_tf(y,y_3)=\sum_{k=k_0}^\infty A_t^kf(y,y_3),
$$
where one finds that 
$$
A^k_t f(y,y_3)=2^{-k|\ka|} \int_{\bR^2} f\Big(y-t\de_{2^{-k}}(x),
y_3-t(1+2^{-k}\phi^k(x)\Big)\,\eta(\de_{2^{-k}}x)\chi(x)\,dx,
$$
with $\phi^k$ is as before. We denote by $\M^k$  the
maximal operator associated to the averaging operators $A^k_t.$ 
Observe next that the  scaling operators $T^k, $  defined by
$$
T^kf(y,y_3):=2^{\frac{k|\ka|}{p}}f(\de_{2^k}(y),y_3),
$$
act isometrically on  $L^p(\bR^3),$  and 
$$
(T^{-k}A^k_tT^k)f(y,y_3)=2^{-k|\ka|} \int_{\bR^2}
f\Big(y-tx,y_3-t(1+2^{-k}\phi^k(x))\Big)\,\eta(\de_{2^{-k}}(x))\chi(x)\,dx,
$$
so that we are reduced to estimating the maximal operator associated to this family of averaging operators. This, in return, can be accomplished by means of Proposition \ref{maxs4.7}, where we choose $\ve:=2^{-k}.$ The resulting estimates can then be summed over $k,$ and we arrive at the desired estimation for $\M$ in this case.

\medskip
Assume next that the coordinates $x$ are {\it not adapted} to $\phi.$ Eventually, we then again arrive at the situation studied in the third step (compare Subsection \ref{near}), i.e., we need to estimate the contribution of the domain $D_\la$ to the maximal operator $\M.$ Again, one can indeed reduce to a small subdomain of the form 
$$
|x_2-\psi(x_1)-cx_1^{a_\la}|\le \ve_0 x_1^{a_\la}, \quad\mbox{ with } x_1>0,
$$
and also assume that \eqref{2.15} holds true at the point $(1,c)=(1,0).$ Recall that we may then write 
$$
\pad_\pr\y=y_2^BQ\y+y_1^n,
$$
where $Q$ is a $\ka^\la$-homogeneous polynomial such that
$Q(1,0)\ne 0,$ and  $B\ge h>2.$
Again, in  this situation, a reduction to the one-dimensional case is no longer possible. 

\medskip
It turns out that this is the most difficult case, which requires various further ideas. The heart of the matter are in fact   rather precise estimations of certain classes of two-dimensional oscillatory integrals depending on small  parameters (cf. \S 5 in \cite{IKM-max}).

\medskip
Indeed, in  order to understand the behavior of $\pad$ as a function of $y_2,$ for $y_1$ fixed, 
we  decompose 
$$
\pad\y=\pad(y_1,0)+\th\y,
$$
and write the complete phase $\Phi$ for $J(\xi)$  in adapted coordinates $\y$ in the form
$$
F(y,\xi)=(\xi_3\pad(y_1,0)+\xi_1 y_1+\xi_2\psi(y_1))\  +\ (\xi_3 \th\y +\xi_2y_2).
$$
Notice that 
$$
\pad(y_1,0)=y_1^n\rho(y_1), \quad \psi(y_1)=y_1^{m_1}\om(y_1),\quad \th_{\ka^\la}\y=y_2^B Q\y,
$$
where $\th_{\ka^\la}$ denotes the $\ka^\la$-principal part  of $\th,$ and where $\rho$ and $\om$ do not vanish at $y_1=0.$

\smallskip
Now, by means of  a dyadic decomposition and re-scaling using the $\ka^\la$-dilations $\{\de^\la_r\}_{r>0}$ we would like to  reduce our considerations as before  to the domain where $y_1\sim 1.$ In this domain, $|y_2|\ll 1,$ so that $\th_{\ka^\la}(y)\sim y_2^BQ(y_1,0).$ What leads to problems is that the ``error term'' $\th_{r}:=\th-\th_{\ka^\la},$ which consists  of terms of higher $\ka^\la$-degree than $\th_{\ka^\la},$ may nevertheless contain terms   of lower $y_2$-degree $l_j<B$ of the form $c_j y_2^{l_j} y_1^{n_j},$ provided $n_j$ is sufficiently large (this corresponds to a ``fine splitting of roots'' of $\theta$ when $\phi$ is analytic). After scaling the $k$-th dyadic piece in our decomposition by $\de^\la_{2^{-k}}$  in order to achieve that $y_1\sim 1$ and $|y_2|\lesssim \ve_0,$ such  terms will have small coefficients compared to the one of $y_2^B Q\y,$ but for $|y_2|$ very small they may nevertheless become dominant and have to be taken into account.

\medskip
In order to resolve this problem, we apply  a {\it further domain decomposition}  by means of a suitable  stopping time argument,  into homogeneous domains $D'_\ell$ and transition domains $E'_\ell,$   oriented, in some sense, at the level sets of 
$\pa_2\pad,$ which in return again are chopped up into dyadic respectively  bi-dyadic pieces. 
 After re-scaling, the contributions of these pieces to the maximal operators  can eventually be estimated by means of  {\it oscillatory integral technics in two variables.} 
 
 \medskip
 More precisely, it turns out that what is needed   are  {\it uniform estimates} for various classes of oscillatory integrals of the form
$$
J(\lambda,\sigma,\delta):=\int_{\bR^2}e^{i\lambda
F(x,\sigma,\delta)}a(x,\delta)\, dx,\qquad (\la>0),
$$
with a  phase function  $F$ of the form 
$$
F(x_1,x_2,\sigma,\delta):=f_1(x_1,\delta)+\si f_2(x_1,x_2,\delta),
$$
and an amplitude $a$ defined for $x$ in some open neighborhood of the origin in $\RR^2$ with compact support in $x.$
The functions $f_1,f_2$  are assumed to be real-valued and depend, like  the function $a,$ smoothly on $x$ and on small real parameters $\de_1,\dots,\de_\nu,$ which form the vector $\de:=(\de_1,\dots,\de_\nu)\in \bR^\nu.$ $\si$ denotes another small real parameter.

With a slight abuse of language we  shall say that $\psi$ is  compactly supported in some open set  $U\subset \RR^2$  if there is a compact subset $K\subset U$ such that $\supp \psi(\cdot,\de)\subset K$ for every $\de.$

To give an idea as to which type of oscillatory integrals we need to estimate, let us remark that the most difficult instance  are ``oscillatory integrals of degenerate Airy type'', for which we have the following result:

\begin{thm}\label{maxs8.3}
Assume that 
$$
|\pa_1f_1(0,0)|+|\pa_1^2f_1(0,0)|+|\pa_1^3f_1(0,0)|\ne 0\ \mbox{ and }\ \pa_1\pa_2f_2(0,0,0)\neq 0,
$$
and that there is some $m\ge 2$ such that 
\begin{equation}\nonumber
\pa_2^lf_2(0,0,0)=0\mbox{  for } l=1,\dots, m-1\mbox{ and  }
\pa_2^{m}f_2(0,0,0)\neq0.
\end{equation}

Then there exist a
neighborhood $U\subset \bR^2$  of the
origin and  constants $\ve,\ve'>0$ such that for any  amplitude $a$ which is compactly supported in $U$   the
following estimate
\begin{equation}\label{max8.9}
|J(\lambda,\sigma,\delta)|\le
\frac{C}{\lambda^{\frac
{1}{2}+\ve}|\si|^{(l_m+c_m\ve)}}
\end{equation}
holds true  uniformly for $|\si|+|\delta|<\ve',$ where $l_m:=\frac{1}{6} $ and $c_m:=1$ for $m<6$, and $l_m:=\frac{m-3}{2 (2m-3)}$ and $c_m:=2$ for $m\geq 6.$ 
\end{thm}

Observe that the order of decay $O(\la^{-1/2-\ve})$ in \eqref{max8.9} is just what we need in order to apply the usual method to control the maximal operator on $L^2$ by means of Sobolev's embedding theorem applied to the scaling variable $t>0.$
\begin{remark}\label{max}
If the Fourier transform of the surface carried measure $\mu$  decays more slowly  than   $O(|\xi|^{-1/2-\ve}),$ then one  cannot directly  apply Sobolev's embedding theorem in order to control the maximal function on $L^2.$ This happens in many situations  when the Gaussian curvature of the hypersurface $S$ vanishes. 
One  method  to overcome this problem is to use suitable damping factors in the amplitude of the corresponding oscillatory integrals, for instance powers of the Gaussian curvature, and combine this with a suitable complex interpolation argument. This technique had been introduced in  \cite{sogge-stein}, and also been used in our first approach our  problems in \cite{ikm}. 

However, the choice of a suitable damping factor can become  quite tricky a task, and we believe that the techniques in \cite{IKM-max}, which avoid damping factors and rather rely on suitable decompositions of the given surface in combination with re-scaling arguments, are simpler and more straight-forward.

For an approach  based on damping methods, we refer to \cite{greenblatt-max1}, where most of those cases are treated by damping techniques which can be reduced to  a one-dimensional problem.
\end{remark}

\setcounter{equation}{0}
\section{Problem C: Fourier restriction to  the hypersurface  $S.$}\label{FREST}

We finally turn to the last problem, namely the determination of  the range of exponents $p$ for which an $L^p$-$L^2$ Fourier restriction estimate 
\begin{equation}\label{rest1}
\Big(\int_S|\hat f(x)|^2\, d\mu(x)\Big)^{1/2}\le C \|f\|_{L^p(\RR^3)}, \qquad f\in \cS(\RR^3),
\end{equation}
holds true.

The idea  of Fourier restriction goes back to Stein,  and a first instance of this concept is the sharp  $L^p$-$L^q$ Fourier restriction estimate  for the circle  in the plane by Zygmund \cite{zygmund}, who extended earlier work by Fefferman and Stein \cite{fefferman} (see also H\"ormander \cite{hoermander}, \cite{carleson-sjoelin} for estimates on more general oscillatory integral operators).  For subvarieties of higher dimension,  the first fundamental result was obtained (in various steps) for  Euclidean spheres $S^{n-1}$ by Stein and Tomas,  who eventually proved that an  $L^p$-$L^2$ Fourier restriction estimate  holds true for $S^{n-1}, n\ge 3,$ if and only if $p'\ge 2(2/(n-1)+1)$ (cf. \cite{stein-book} for the history of this result).

\medskip
Even though I shall not pursue that problem here, let me  briefly remind that  a more general and  even  substantially deeper problem is to determine the exact range of exponents $p$ and $q$ for which an $L^p$-$L^q$ restriction estimate 
$$
\Big(\int_{S^{n-1}}|\hat f(x)|^q\, d\si(x)\Big)^{1/q}\le C \|f\|_{L^p(\RR^n)}
$$
holds true for spheres. It is conjectured that this is the case  if and only if $p'>2n/(n-1)$ and $p'\ge q(2/(n-1)+1), $  and there has been  a lot of very deep work on this problem by many  mathematicians, including  Bourgain, Wolff,  Vargas, Vega, Katz, Tao, Keel, Lee, and most recently Bourgain and Guth \cite{bourgain-guth},   which has led to important  progress, but the problem is still open in dimensions $n\ge 3$ and represents one of the major challenges in Euclidean harmonic analysis, bearing various deep connections with  other  important open problems, such as the  Bochner-Riesz conjecture, the Kakeya conjecture and  Sogge's local smoothing conjecture for solutions to the wave equation. I refer to  Stein's book  \cite{stein-book} for more information on these topics and their history until 1993, and various related essays by T. Tao  which  can be found on his webpage.
\medskip

Coming back to the restriction estimate  \eqref{rest1} for our hypersurface $S$ in $\RR^3,$ we begin with  the case where  there exists a linear coordinate system which is adapted to the function $\phi.$ For this case, a  complete answer had been given in   \cite{IM-uniform} (for analytic hypersurfaces, partial results had been obtained before by Magyar \cite{magyar}) :
\begin{thm}\label{adarestrict}
Assume that, after applying a suitable linear change of coordinates, the coordinates $\x$ are adapted to $\phi.$ 
We then define the critical exponent $p_c$ by
\begin{equation}\label{pcritical}
p'_c:=2h(\phi) +2,
\end{equation}
where $p'$ denotes the exponent conjugate to $p,$ i.e., $1/p+1/p'=1.$

Then there exists a neighborhood $U\subset S $ of the point $x^0$ such
that for every  non-negative density $\rho\in C_0^\infty(U)$ the  Fourier restriction estimate \eqref{rest1}
holds true for every $p$ such that
\begin{equation}\label{rest2}
1\le p\le p_c.
\end{equation}

Moreover, if $\rho(x^0)\ne 0,$ then the condition \eqref{rest2} on $p$ is also necessary for the validity of  \eqref{rest1}.
\end{thm}

In many cases, this result is an immediate consequence of Theorem \ref{estu3} and Greenleaf's  classical restriction estimate in \cite{greenleaf}, which I shall state here  for the special case of hypersurfaces  only.
 \begin{thm}[Greenleaf]\label{greenleaf}
Assume that $|\widehat{d\mu}(\xi)|\lesssim |\xi|^{-1/h}.$ Then the restriction estimate \eqref{rest1} holds true for every $p\ge 1$ such that $p'\ge 2(h+1).$
\end{thm}
Indeed, observe that a problem arises only when Varchenko's exponent $\nu(\phi)$ equals one  in \eqref{estu2}. In this case, a direct application of Greenleaf's result yields only  the range $1\le p<p_c.$ 

To capture also the endpoint $p=p_c,$ recall from Subsection \ref{adaptedc} that we had effectively  decomposed the measure $\mu$ into a dyadic sum $\mu=\sum_{k=k_0}^\infty \mu_k,$  where $\mu_k=(\chi_k\otimes 1)\mu,$ and where $\hat\mu_k(\xi)=J_k(\xi)$ is given by \eqref{intk1}. Moreover, we had estimated $J_k(\xi)$ in \eqref{estjk}, which in particular yields the estimate
$$
 |\widehat{\mu_k}(\xi)|\le C 2^{-k|\ka|}(1+2^{-k}|\xi_3|)^{-1/h}.
$$
The measures $\mu_k$ are supported in dyadic annuli of ``radius'' $2^{-k},$ which are images of an annulus of radius of size one  under the dilations $\de_{2^{-k}}$ associated to the principal weight $\ka,$ so that we cannot directly apply Greenleaf's restriction estimate. However, it is important to notice that our estimates for 
$\widehat{\mu_k}(\xi)$ do not carry  a logarithmic factor yet (that only arose for $\hat \mu(\xi)$ in certain cases through summation over the  $k$), and a simple re-scaling argument can then be applied to derive the following uniform restriction estimate for the family of measures $\mu_k:$
\begin{eqnarray}\label{rest5.4}
\int|\hat f(x)|^2\,d \mu_{k}(x)\le C^2 \|f \|_{p_c}^2\, .
\end{eqnarray}

Fix a cut-off function $\tilde\chi\inÊC^\infty_0(\RR^2)$ supported in an annulus centered at the origin
such that $\tilde\chi=1$ on the support of $\chi,$ and define dyadic frequency decomposition operators $\Delta'_k$ by
\begin{eqnarray*}
\widehat{\Delta'_kf}(x):=\tilde\chi(\de_{2^k}x')\, \hat f(x',x_3),
\end{eqnarray*}
where we have written $(x_1,x_2)=x'.$
 Then $\int|\hat f (x)|^2d\mu_k(x)=\int |\widehat{\Delta'_kf}(x)|^2d\mu_k(x),$  and setting $p:=p_c,$ we see that \eqref{rest5.4} yields in fact 
\begin{eqnarray*}
\int|\hat f(x)|^2d\mu_k(x)\le C^2\, \|\widehat{\Delta'_kf}\|_{p}^2,
\end{eqnarray*}
for any $k\ge k_0$.  In combination with Minkowski's inequality, this implies
\begin{eqnarray*}
\left(\int|\hat f(x)|^2d\mu(x)\right)^{1/2}=\left(\sum_{k\ge k_0}\int|\hat f(x)|^2d\mu_k(x)\right)^{1/2}\le \left(\sum_{k\ge k_0}\|\Delta'_kf\|_{p}^2\right)^{1/2}\\
=C\left(\left(\sum_{k\ge k_0}\left(\int|\Delta'_kf(x)|^{p}dx\right)^{2/p}\right)^{p/2}\right)^{1/p}\le
C\left\|\left(\sum_{k\ge k_0}|\Delta'_kf(x)|^2\right)^{1/2}\right\|_{L^{p}(\bR^3)},
\end{eqnarray*}
since $p<2$.
Estimate \eqref{rest1} for $p=p_c$ follows thus by means of Littlewood-Paley theory.

\bigskip

From now on, we shall therefore always make the following 
\begin{assumption}\label{assumption}
There is no linear coordinate system which is adapted to $\phi.$
\end{assumption}

According to our discussion in Subsection \ref{consada}, we may then also assume that the coordinates $x$ are linearly adapted to $\phi,$ and  that there are adapted coordinates $y$ of the form $y_1=x_1, y_2=x_2-\psi(x_1),$ where 
\begin{equation}\label{rest5.5}
\psi(x_1)=x_1^{m}\om(x_1),\quad \mbox{ with } \  \om(0)\ne 0 \ \mbox{Êand  } \ m\ge 2.
\end{equation}
$\pad(y)$ will again denote $\phi$ when expressed in these adapted coordinates, and we shall use the notions introduced for the study of the Newton polyhedron of $\pad$ in Subsection \ref{nonadaptedc}.

\medskip

Consider the line parallel to the bi-sectrix
$$
\Delta^{(m)}:=\{(t,t+m+1):t\in\RR\}.
$$
For any edge $\ga_l\subset L_l:=\{(t_1,t_2)\in \bR^2:\ka^l_1t_1+\ka^l_2 t_2=1\}$  define $h_l$ by
$$
\Delta^{(m)}\cap L_l=\{(h_l-m, h_l+1)\},
$$ i.e.,
\begin{equation}\label{hl}
h_l=\frac {1+m\ka^l_1-\ka^l_2}{\ka^l_1+\ka^l_2},
\end{equation}
and define the {\it restriction height}, or short, {\it $r$-height,} of $\phi$ 
\color{black} by 
$$
h^r(\phi):= \max(d, \max\limits_{\{l=1,\dots, n+1:a_l>m\}} h_l).
$$

  %%%%%%%%%%%%%%%%%%%%%%
%%%%%%%%%%%%%%%%%%%%%%
%%%%%%% Figure 5 %%%%%%%%%
%%%%%%%%%%%%%%%%%%%%%%
%%%%%%%%%%%%%%%%%%%%%%
\begin{figure}
\centering
\scalebox{0.35}{\input{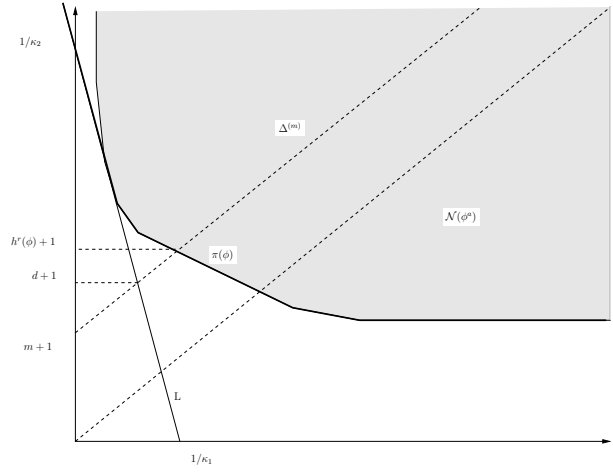}}
\caption{r-height}{\label{fig5}}
\end{figure}

\begin{remarks}\label{r1}
 \begin{itemize}
\item[(a)] For $L$ in place of $L_l$ and $\ka$ in place of $\ka^l,$  one has $m=\ka_2/\ka_1$ and $d=1/(\ka_1+\ka_2),$ so that one gets $d$ in place of $h_l$ in \eqref{hl}.

\item[(b)]  Since $m< a_l,$ we have $h_l<1/(\ka^l_1+\ka^l_2),$ hence $h^r(\phi)<h(\phi).$
\end{itemize}
\end{remarks}

It is easy to see by Remark \ref{r1} (a)  that the  $r$-height admits the following {\it geometric interpretation:}
\medskip

By following Varchenko's algorithm (cf. Subsection 8.2 of \cite{IKM-max}), one realizes that the principal line $L$ of $\N(\phi)$ is a supporting line also for the Newton polyhedron of $\pad,$ which intersects  $\N(\pad)$ in a compact face,  either in a single vertex, or a compact edge. I.e., the intersection contains at least one and at most two vertices of $\N(\pad),$ and we choose  $(A_{l_0-1},B_{l_0-1})$ as the one with smallest second coordinate.  Then $l_0$ is the smallest index $l$ such that $\ga_l$ has a slope smaller than the slope of $L,$  i.e., 
 $
a_{l_0-1}\le m<a_{l_0}.
$
We may thus consider  the {\it augmented} Newton polyhedron $\N^r(\pad)$ of $\pad,$ which is the convex hull of the union of $\N(\pad)$ with  the half-line $L^+\subset L$ with right endpoint $(A_{l_0-1},B_{l_0-1}).$ Then $h^r(\phi)+1$ is the second coordinate of the  point at which the line $\Delta^{(m)}$ intersects the boundary of $\N^r(\pad).$ 
\medskip

The main result from \cite{IM-rest1}, \cite{IM-rest2} then reads as follows.

 \begin{thm}\label{nonadarestrictanal}
Let $\phi\ne 0$ be real analytic, and assume that there is no linear  coordinate system  adapted to $\phi.$  Then there exists a neighborhood $U\subset S $ of  $x^0=0$ such
that for every  non-negative density $\rho\in C_0^\infty(U),$ the Fourier restriction estimate \eqref{rest1} holds true for every $p\ge 1$ such that 
$p'\ge p'_c:=2h^r(\phi)+2.$
\end{thm}

\begin{remarks}\label{r2}
 \begin{itemize}
\item[(a)] An application of Greenleaf's result would  imply, at best,  that the condition  $p'\ge  2h(\phi)+2$ is sufficient for \eqref{rest1} to hold, which is a strictly stronger condition than  $p'\ge  p_c'.$  
\item[(b)] A. Seeger recently informed me that in a preprint, which  regretfully had remained unpublished,   Schulz \cite{schulz-un} had already observed this  kind of phenomenon for particular examples  of surfaces of revolution.
\item[(c)] It can be shown that the number $m$ is well-defined, i.e., it does not depend on the chosen linearly adapted  coordinate system $x.$
\color{black}
 \end{itemize}
\end{remarks}

\color{black}

\begin{example}\label{ex1}
{\rm $$\phi\x:=(x_2-x_1^{m})^n, \qquad n,m\ge 2.$$
The coordinates $\x$ are not adapted. Adapted coordinates are 
$y_1:=x_1, y_2:=x_2-x_1^m,$ in which $\phi$ is given by
 $$\phi^a(y_1,y_2)=y_2^n.$$
Here 
 \begin{eqnarray*}
&&\ka_1=\frac 1{mn}, \quad \ka_2=\frac 1n,\\
&&d:=d(\phi)= \frac 1{\ka_1+\ka_2}=\frac {nm}{m+1}<n,
\end{eqnarray*}
and 
$$ p'_c=\left\{  \begin{array}{cc}
2d+2, & \mbox{   if } n\le m+1 ,\hfill\\ 
2n, & \mbox{   if } n>m+1\ .\hfill 
\end{array}\right.
$$
On the other hand, $h:=h(\phi)=n,$ so that $2h+2=2n+2>p'_c.$ }
\end{example}

An analogous theorem holds true even for smooth, finite type functions $\phi,$ under an  additional {\it Condition (R)} which, roughly speaking, requires that whenever  the Newton diagram suggests that a root with leading term given by the principal root jet $\psi(x_1)$ should have multiplicity $B,$ then indeed such a root of multiplicity $B$ 
does exist (this is a condition on the behavior of flat terms). Condition (R)  is always satisfied when $\phi$ is real-analytic.

For example, Condition (R) would hold true  for 
$$
\phi_g\x= (x_2-x_1^2-f(x_1))^2,
$$
for every flat  smooth function $f(x_1)$ (i.e., $f^{(j)}(0)=0$ for every $j\in\NN$).
On the other hand, (R) is not satisfied for 
$$
\phi_b\x:= (x_2-x_1^2)^2+ f(x_1),
$$
unless $f$ vanishes identically.

\medskip
I  also like to mention that there is a more invariant description of the notion of $r$-height (cf. Proposition 1.9 in \cite{IM-rest1}),   somewhat in the spirit of Varchenko's definition of height, but I refrain from stating it here since this would require  the introduction of further, somewhat technical notions.

\subsection{Necessity of the condition $p'\ge 2h^r(\phi)+2$ }\label{knapp}

In order to better understand the meaning of the notion of r-height, let me present the proof of the necessity of  the condition $p'\ge 2h^r(\phi)+2$  for the validity of the Fourier restriction estimate \eqref{rest1} when $\rho(x^0)\ne 0.$ The proof will be based on a  modified Knapp-type argument. 

\medskip
Let $\ga_l$ be any edge of $\N(\pad)$ with $a_l>m,$ and choose the weight $\ka^l$ such that $\ga_l$ lies on the line $L_l$ given by $\ka^l_1t_1+\ka^l_2t_2=1.$ Consider the region
$$D^a_\ve:=\{y\in\RR^2:|y_1|\le\ve^{\ka^l_1}, |y_2|\le\ve^{\ka^l_2}\},\quad\ve>0,
$$ 
in adapted coordinates $y.$ In the original coordinates $x$, it corresponds to
$$
D_\ve:=\{x\in\RR^2:|x_1|\le\ve^{\ka^l_1}, |x_2-\psi(x_1)|\le\ve^{\ka^l_2}\}.
$$
Assume that $\ve$ is sufficiently small. Since 
$$
\pad(\ve^{\ka^l_1}y_1,\ve^{\ka^l_2}y_2)=\ve \,\Big(\pad_{\ka^l}(y_1,y_2)+O(\ve^\delta)\Big)
$$
for some $\de>0,$ we have that $|\pad(y)|\le C\ve$ for every $y\in D^a_\ve,$ i.e.,
\begin{equation}\label{n1}
|\phi(x)|\le C \ve \quad \mbox {   for every   }x\in D_\ve.
\end{equation}
Moreover, for $x\in D_\ve,$
$$
|x_2|\le \ve^{\ka^l_2}+|\psi(x_1)|\lesssim \ve^{\ka^l_2}+\ve^{m\ka^l_1}.
$$

Since $m\le a_l=\ka^l_2/\ka^l_1,$ we find that
$$
|x_2|\lesssim \ve^{m\ka^l_1},
$$
so that we may assume that  $D_\ve$ is contained in the box where $|x_1|\le\ve^{\ka^l_1}, |x_2|\le\ve^{m\ka^l_1}.$ Choose $f_\ve$ such that
$$
\widehat{f_\ve}(x_1,x_2,x_3)=\chi_0\Big(\frac {x_1}{\ve^{\ka^l_1}}\Big)\chi_0\Big(\frac {x_2}{\ve^{m\ka^l_1}}\Big)\chi_0\Big(\frac {x_3}\ve\Big).
$$
Then by \eqref{n1} we see that $\widehat{f_\ve}(x_1,x_2,\phi(x_1,x_2))\ge 1$ on $D_\ve,$ hence, if $\rho(0)\ne 0,$ then 
$$
\Big(\int_S |\widehat {f_\ve}|^2\, \rho d\si\Big)^{1/2}\ge |D_\ve|^{1/2}=\ve^{(\ka^l_1+\ka^l_2)/2}.
$$
Since $\|f_\ve\|_p\simeq \ve^{((1+m)\ka^l_1+1)/p'},$ we find that the restriction estimate can hold true only if 
$$
p'\ge 2\frac{(1+m)\ka^l_1+1}{\ka^l_1+\ka^l_2}=2h_l+2,
$$
where we recall that $h_l=(1+m\ka^l_1-\ka^l_2)/(\ka^l_1+\ka^l_2).$

Notice that the argument still works if we replace the previous line $L_l$  by the line $L$ associated to the weight $\ka,$ and $\pad_{\ka^l}$ by $\pad_{\ka}.$ Since here $m\ka_1=\ka_2,$ this  leads to the condition $p'\ge 2d+2,$ so that altogether 
necessarily 
$$
p'\ge 2 \max (d,\max\limits_{l: a^l>m} h_l) +2=2 h^{r}(\phi)+2.
$$
\qed

\subsection{Sufficiency of the condition $p'\ge2h^{r}(\phi)+2$: I. Some key steps in the proof when $\hl(\phi)\ge 2$ }\label{hlinge2}

Let us assume that we are working in linearly adapted coordinates $x,$ so that $d:=d(\phi)=\hl:=\hl(\phi).$ 

In the preceding discussions of problems A and B, it had been natural to distinguish between the cases where $h:=h(\phi)<2$ and where $h\ge 2,$ since in the latter case, in many situations a reduction to a one-dimensional situation had been possible by means of the van der Corput type Lemma \ref{corput}. For similar reasons, in the discussion of  Problem C it appears natural to distinguish between the cases where $d<2$ and where $d\ge 2.$ In addition, when $d\ge 2,$  it turns out that  the case where $d\ge 5$ can be handled in a somewhat simplified way compared to the case  where $2\le d<5.$ 

\medskip
We shall therefore assume in this section that  $d\ge 5.$  In a {\it first step}, we can again localize to the  narrow   $\ka$-homogeneous subdomain \eqref {restdomain} of the curve $x_2=b_1x_1^m$ given by 
$$
|x_2-cx_1^m|\le \ve x_1^m,
$$
by means of a cut-off function $\rho_1.$ Indeed, the technique of proof that we used in the case of adapted coordinates can essentially be carried  over to the domain complementary to  \eqref {restdomain} without major new ideas, since one can show that  the Fourier transforms of the corresponding dyadic pieces $\mu_k$ of the measure $\mu$ satisfy estimates of the form
$$
 |\widehat{\mu_k}(\xi)|\le C 2^{-k|\ka|}(1+2^{-k}|\xi_3|)^{-1/d}.
$$
Notice also that $h^r:=h^r(\phi)\ge d.$ 

\medskip
Let us again assume for instance that  the principal face of the Newton polyhedron of $\pad$ is a compact edge. Using the same notation as in Section \ref{FT}, we choose  again $\la>l_0$ so  that the edge $\ga_\la=[(A_{\la-1},B_{\la-1}) ,(A_\la,B_\la)]$ is the principal face $\pi(\pad)$ of the Newton polyhedron of $\pad.$

\medskip
In a {\it second step}, we again narrow down the domain \eqref{restdomain} to the  neighborhood 
$D_\pr:=D_\la$ of  the principal root jet given by \eqref{restdomain2},  where
$$
 |x_2-\psi(x_1)|\le N_\la  x_1^{a_\la},
$$
again by decomposing the difference set of the  domains \eqref{restdomain} and \eqref{restdomain2}  into   the domains
 $$
D_l:=\{\x:\ve_l  x_1^{a_l}< |x_2-\psi(x_1)|\le N_l x_1^{a_l}\},\quad l=l_0,\dots,\la-1,
$$
 and  the intermediate  domains
$$E_l:=\{\x:N_{l+1} x_1^{a_{l+1}}<|x_2-\psi(x_1)| \le \ve_l x_1^{a_l}\}, \quad l=l_0,\dots,\la-1,
$$
as well as 
$E_{l_0-1}:=\{\x:N_{l_0} x_1^{a_{l_0}}<|x_2-\psi(x_1)| \le \ve_1 x_1^{m}\}.
$
\medskip

\noi{\bf Contribution by the domains $E_l.$}
Let us denote by $\mu_{E_l}$ the contribution of the  transition domains $E_l$  to the measure $\mu.$ We then decompose $\mu_{E_l}$ bi-dyadically w.r. to the adapted coordinates $y$ as 
$$
\mu_{E_l}=\sum_{j,k}\mu_{j,k},
$$
so that $\mu_{j,k}$ is supported where $y_1=x_1\sim 2^{-j}$ and $y_2=x_2-\psi(x_1)\sim 2^{-k}.$ Observe that this a  ``curved rectangle'' in our original coordinates $x.$ In a similar way as in the case of adapted coordinates, we would like to localize to these curved rectangles by means of Littlewood-Paley theory in order to reduce to uniform restriction estimates for the family of measure $\mu_{j,k},$ i.e., 
\begin{equation}\label{6.5}
\int_S |\widehat f|^2\,d\mu_{j,k}\le C\|f\|^2_{L^{p}(\RR^3)},
\end{equation}
for $p\le p_c.$  Clearly, because of the non-linearity $\psi(x_1),$ this is not possible by means of Littlewood-Paley techniques in the variables $x_1$  and $x_2,$ but it turns out the we can use the variables $x_1$ and $x_3$ to accomplish this.

Indeed, one can show that 
$$
\pad(y)=c_l\,y_1^{A_l}y_2^{B_l}\Big(1+\mbox{small error}\Big) \quad \mbox{on}\quad E^a_l,
$$
which in return implies that on the domains $E_l$ respectively $E^a_l$ (recall that  $E^a_l$ represents $E_l$ in the adapted coordinates $y$) the conditions $y_1\sim 2^{-j}, y_2\sim 2^{-k}$ are equivalent to the conditions
$$
x_1\sim 2^{-j} \quad \mbox{and} \quad \phi(x)\sim 2^{-(A_lj+B_lk)}
$$
(cf. Lemma 6.1 in \cite{IM-rest1}).

Working in the coordinates $y,$ after re-scaling of  the measures  $\mu_{j,k}$ to get normalized measures $\nu_{j,k}$ supported on a surface $S_{j,k}$ where $y_1\sim 1\sim y_2,$ by means of the formula above one eventually finds that $S_{j,k}$ is a 
 small perturbation of the limiting surface
$$
S_{\infty}:=\{(y_1, \, y_1^m\om(0),\, cy_1^{A_l}y_2^{B_l}):  y_1\sim 1\sim y_2\},
$$
But  $|\pa(cy_1^{A_l}y_2^{B_l})/\pa y_2|\sim 1,$ since $B_l\ge 1,$  which shows that $S_{\infty},$ and hence also $S_{j,k},$ is a 
smooth hypersurface with one non-vanishing principal  curvature (with respect to $y_1$) of size $\sim 1.$  This implies that
$$
|\widehat{ \nu_{j,k}}(\xi)|\le C (1+|\xi|)^{-1/2},
$$
uniformly in $j$ and $k.$ Applying Greenleaf's restriction theorem to these measures, and scaling these estimates back, we eventually arrive (in a not completely trivial way) at the estimates \eqref{6.5}. 
It is important to observe here that Greenleaf's results implies restriction theorems for $p'\ge 2(1+2)=6,$ which is sufficient for our purposes, since $p'_c\ge 2d+2,$ where $d\ge 2.$

\medskip
\noi{\bf Contribution by the domains $D_l.$}
Let us next turn to the domains $D_l.$ After dyadic decomposition of the domain $D_l$ in the adapted coordinates $y$ by means of the $\ka^l$-dilations and suitable re-scaling, the re-scaled measure $\nu_k$ corresponding to the measures $\mu_k$ turns out to be of the form
\begin{eqnarray*}
\laa\nu_k,\, f\ra :=\int f(y_1,\, 2^{(m\ka_1^l-\ka_2^l) k}y_2+y_1^m\om(2^{-\ka_1^lk} y_1),\, \phi^k(y))
\, \tilde\eta(y)\, dy,
\end{eqnarray*}
and by means of a finite partition of unity, we may assume that the amplitude $\tilde \eta$ is supported in a sufficiently thin set $U(c_0),$ on which
$$
y_1\sim 1\quad\mbox{and}\quad |y_2-c_0y_1^{a_l}|\le \ve y_1^{a_l}.
$$
This measure $\nu_k$ is supported in a variety $S_k$ which in the limit as $k\to\infty$ tends to the variety 
$$
S_\infty:=
\{g_\infty\y:=(y_1,\,\om(0)y_1^m\,, \pad_{\ka^l}(y)): \y\in U(c_0)  \},
$$
since $m\ka_1^l-\ka_2^l<a_l\ka_1^l-\ka_2^l=0$ and since $\phi^k$ tends to $\pad_{\ka^l}.$ Here, $c_0$ is fixed with $|c_0|\le N_l.$

Again, we have to prove uniform restriction estimates for the family of measures $\nu_k.$ Depending on $c_0,$ different cases may arise.

\medskip
{\bf 1. Case.} $\pa_2\pad_{\ka^l}(1,c_0) \neq 0$.  Then we may use $z_2:=\pad_{\ka^l}(y_1,\, y_2)$ in place of $y_2$ as a new coordinate for $S_\infty$ (which thus is a hypersurface), and  since $y_1\sim 1$ on $U(c_0),$ we find  that $S_\infty,$  hence also $S_k,$ is a hypersurface with {\it one non-vanishing principal curvature.} Then we may essentially argue as for the domains $E_l.$

\medskip
{\bf 2. Case.} $\pa_2\pad_{\ka^l}(1,c_0)=0,$ but $\pa_1\pad_{\ka^l}(1,c_0)\neq 0$. In this case, since $\pad_{\ka^l}$ is a $\ka^l$-homogenous polynomial, by Euler's homogeneity relation we have also $\pad_{\ka^l}(1,\, c_0)\neq0$. One can then show that one can fibre the variety $S_\infty$ into the family of curves 
$$
\ga_c(y_1) := g_\infty(y_1,cy_1^{a_l})=(y_1,\om(0)y_1^m,\pad_{\ka^l}(y_1,cy_1^{a_l})),
$$
for $c$ sufficiently close to $c_0,$ and one finds that the curve $\ga_{c_0}(y_1)=(y_1,\om(0)y_1^m,b_0y_1^{1/\ka_1^l})$ has  {\it non-vanishing torsion,} since $b_0\ne 0.$ The same applies then to  the curves $\ga_c,$  and for $k$ sufficiently large, we do obtain the analogous results for the varieties $S_k.$ 

This allows to decompose the measure $d\nu_k$ as a direct integral of measures $d\Gamma_c$ supported on curves $\ga^l_c$ with non-vanishing torsion.

We may  thus apply Drury's Fourier restriction theorem for curves with non-vanishing torsion (cf. Theorem 2 in \cite{drury} and  \cite{bak-oberlin-seeger}, \cite{dendrinos-m}) to the measures
$d\Gamma_c$ and obtain uniform estimates 
$$
\Big(\int |\hat f|^2\, d\Gamma_c\Big)^{\frac 12}\le C_{p}\|f\|_{L^p(\RR^3)},
$$
when $p'>7$ and $2\le p'/6.$ Since we assume here that $p_c'\ge 2(d+1)>2(5+1)=12,$ these estimates, after re-scaling to the measures $\mu_k,$ yield the desired restriction estimates for the contributions by the domains $D_l.$ 

Notice that it is here that we need the condition $d=\hl>5.$ 

\medskip

{\bf 3. Case.} $\pa_2\pad_{\ka^l}(1,c_0)=0$ and $\pa_1\pad_{\ka^l}(1,c_0)=0$. Then Euler's homogeneity relation implies that also $\pad_{\ka^l}(1,c_0)=0,$ so that $\pad_{\ka^l}$ has a real root of multiplicity $B\ge 2$ at $(1,c_0),$ and one finds that 
\begin{equation}\label{7.6}
\pad_{\ka^l}(y_1,\, y_2)=y_2^{B_l}(y_2-c_0y_1^{a_l})^B Q\y,
\end{equation}
where $Q$ is a $\ka^l$-homogenous  smooth function   such that $Q(1,c_0)\neq0$ and $Q(1,0)\ne 0.$ One can also prove that $B<d/2.$ 
\medskip

We can then essentially follow the Stein-Tomas method for proving  $L^p$-$L^2$- restriction estimates. We localize to frequencies  of size $\La> 1$ by putting 
$$
\widehat{\nu_k^\La}(\xi):=\chi_1\Big(\frac \xi{\La}\Big)\widehat{\nu_k}(\xi),
$$
 where $\chi_1$ is  a smooth bump function supported where $|\xi|\sim 1.$ We claim that the measures $\nu_k^\La$ satisfy the following estimates, uniformly in $k\ge k_0,$ provided $k_0$ is sufficiently large and $\ve'$ sufficiently small:
\begin{eqnarray}\label{7.8}
\|\widehat{\nu_k^\La}\|_\infty&\le& C \La^{-1/B}\,;\\ 
\|\nu_k^\La\|_\infty&\le& C \La^{2-1/B}\,. \label{7.9}
\end{eqnarray}
Indeed,  
$$
\widehat{\nu_k^\La}(\xi)=\chi_1\Big(\frac \xi{\La}\Big)\, \int  e^{-i\Big[\xi_1y_1+\xi_2 \Big(2^{(m\ka_1^l-\ka_2^l) k}y_2+y_1^m\om(2^{-\ka_1^lk} y_1)\Big)+ \xi_3\phi_k(y)\Big]}
\, \tilde\eta(y)\, dy,
$$
which, in the limit as $k\to\infty,$ simplifies as 
$$
\widehat{\nu_\infty^\La}(\xi)=\chi_1\Big(\frac \xi{\La}\Big)\,\int  e^{-i[\xi_1y_1+\xi_2\om(0)y_1^m+ \xi_3\pad_{\ka^l}(y)]}
\, \tilde\eta(y)\,  dy.
$$
Now, if $|\xi_3|\ge c |(\xi_1,\xi_2)|,$ then an application of van der Corput's lemma to the integration in $y_2$ yields 
$|\widehat{\nu_\infty^\La}(\xi)|\lesssim |\xi_3|^{-1/B}$ (cf. \eqref{7.6}), and if $|\xi_3|\ll |(\xi_1,\xi_2)|,$ we may apply van der Corput's lemma to the $y_1$-integration and obtain $|\widehat{\nu_\infty^\La}(\xi)|\lesssim |(\xi_1,\xi_2)|^{-1/2}.$ Since $B\ge 2,$ and because  van der Corput's estimates are  stable under small perturbations, we thus obtain \eqref{7.8}.

In order to verify  \eqref{7.9}, observe that $\nu_\infty^\La(x_1,x_2,x_3)$ is given by 
\begin{eqnarray*}
\La^3\int (\F^{-1}\chi_1)(\La(x_1-y_1),\La(x_2-\om(0)y_1^m), \La(x_3-\pad_{\ka^l}(y_1,y_2))\tilde\eta(y)\, dy_1dy_2,
\end{eqnarray*}
hence, by a change of coordinates,
\begin{eqnarray*}
|\nu_\infty^\La(x_1,x_2,x_3)|
&\le&\La^2\int \rho(z_1)\,\rho( \La(x_3-\pad_{\ka^l}(x_1-\frac{z_1}\La,y_2))\,\eta_1(x_1-\frac{z_1}\La,y_2)\, dz_1dy_2, 
\end{eqnarray*}
where $\rho$ and $\eta_1$ are  suitable, non-negative Schwartz functions, and $\eta_1$ localizes again to $U(c_0).$ However, since $|\pa_2^B\pad_{\ka^l}(y_1,y_2))|\simeq 1$ on the domain of integration, classical  sublevel estimates, originating in work by van der Corput \cite{vdC} (see also \cite{arhipov}, and \cite{ccw},\cite{grafakos}),  essentially imply that the integral with respect to $y_2$ can be estimated by $O(\La^{-1/B}),$ uniformly in $y_1$ and $\La.$

Interpolating  the estimates \eqref{7.8} and \eqref{7.9}, and applying the standard Stein-Tomas argument (see, for instance, \cite{greenleaf}), it is easily seen that we can sum the corresponding estimates over all dyadic $\La\gg1,$ and we obtain the  $L^p$-$L^2$ restriction estimate
$$
\Big(\int |\widehat f|^2\,d\nu_k\Big)^{1/2}\le C_{p} \|f\|_{L^p}
$$
whenever $p'> 4B,$ uniformly in $k,$ for $k$ sufficiently large.  Since $B<d/2,$ we have $ p'_c\ge 2d+2>4B,$ so that the range $p>4B$ does include the critical value $p=p_c.$

\medskip
%\color{red}
Scaling back to the measures $\mu_k,$ we find  that also the original measures $\mu_k$ satisfy a uniform restriction estimate
$$
\Big(\int |\widehat f|^2\,d\mu_k\Big)^{1/2}\le C_{p} \|f\|_{L^p},
$$
where $C_p$ does not depend on $k,$ provided $p'\ge 2h_l+2.$ However, this applies to $p_c,$ since $h^r(\phi)\ge h_l.$

Finally, observe that we can achieve our dyadic decomposition into the measures $\mu_k$ by means of a dyadic decomposition in the variable   $x_1,$ so that these uniform estimates allow to sum over all $k$ by means of Littlewood-Paley theory applied to variable $x_1!$

%\color{black}

\bigskip What remains to be understood is the contribution by the domain
$D_\pr=D_\la$ given by 
$$
 |x_2-\psi(x_1)|\le N_\la  x_1^{a_\la}.
$$
In this domain, the upper bound $B<d/2$ for the multiplicity $B$ of real roots will in general no longer be valid, as  examples show,  not even the weaker condition $B<h^r(\phi)/2,$ which would still suffice for the previous argument.

\medskip
In order to resolve this problem, we  apply  again a {\it further domain decomposition}  by means of a stopping time argument. 
 
Notice first that  if we are to proceed as for the  case $l<\la,$ a major problem arises only in Case 3 where  $\nabla\phi^a_{\pr}(1,c_0)=0;$ in all other cases we can essentially  argue as before.

\smallskip
We therefore devise the {\it stopping  time argument} essentially as follows:
\smallskip

We put $\phi^{(1)}:=\pad.$ 
If Case 3 does not appear for any choice of $c_0,$ then we stop our algorithm  with $\phi^{(1)},$ and are done. Otherwise, if Case 3 applies to $c_0,$ so that  $c_0 y_1^{a_\la}$ is a root of $\pad_{\ka^\la},$    then we  define new coordinates $z$ in place of $y$ by putting
\begin{equation}\nonumber
z_1:=x_1\qquad \mbox{and  }\ z_2:=x_2-\psi(x_1)-c_0 x_1^{a_\la},
\end{equation}
and express $\phi$ by $\phi^{(2)}$ in the coordinates $z.$  Again, if Case 3 does not appear  (for $\phi^{(2)}$ in  place of $\phi^{(1)}$) in the corresponding $z$-domain, we stop our algorithm. We also stop the algorithm  when no further fine splitting of roots does occur, i.e., when there is a root with leading term $\psi(x_1)+c_0 x_1^{a_\la}$ of $\phi$ which has the same  multiplicity as the  trivial root $z_2=0$  of the principal part of $\phi^{(2)}.$   Otherwise, we continue in an analogous way.

This algorithm will stop after a finite number of steps, and  eventually leads to a further domain decomposition  of $D_\pr$  into ``homogeneous'' domains $D_{(l)}$ and transition domains $E_{(l)},$   which  can eventually be treated by methods similar to those applied for the domains $E_l$ and $D_l.$

\medskip
\subsection{Sufficiency of the condition $p'\ge2h^{r}(\phi)+2$: II. Brief sketch of some ideas of the proof when $\hl(\phi)<2$}\label{hlinge3}

This case turns out to be by far  more  difficult than the case where $\hl(\phi)\ge 5,$  and its discussion requires  numerous further methods and ideas which I can  sketch only very briefly.

A first observation is that   $h^r(\phi)=d$  when $d=\hl(\phi)<2,$ so that
$$
p_c'=2d+2.
$$

\medskip

The starting point of our analysis is the following local normal form of our function $\phi.$ It is closely related to the classification of singularities in   \cite{arnold} and \cite{duistermaat} .%and \cite{sirsma}.

\begin{thm}\label{normform}
If $\hl(\phi)< 2,$ then  locally $\phi$ is of the form
\begin{equation}\label{A}
\phi(x_1,x_2)=b(x_1,x_2)(x_2-\psi(x_1))^2 +b_0(x_1).
\end{equation}
Here $b,b_0$ and $\psi$ are  smooth, and $\psi$ is again the principal root jet,  and either
 \begin{itemize}
\item[(a)]  $b(0,0)\ne 0,$ and either $b_0$ is flat (singularity of type $A_\infty$), or of finite type $n,$ i.e.,
$b_0(x_1)=x_1^n\beta(x_1),$  where $\beta(0)\ne 0$ (singularity of type $A_{n-1}$); 
\item[] or 
\item[(b)] $b(0,0)= 0$ and 
$b\x=x_1 b_1\x+x_2^2 b_2(x_2),$ with $b_1(0,0)\ne 0$ \newline(singularity of type $D$). 
 \end{itemize}
\end{thm}

Let us just consider the case where $\phi$ is of finite type $A_{n-1}$ (type D can be treated in a rather similar way). 

In a {\it first step}, by making use of these normal forms in order to estimate certain two-dimensional oscillatory integrals  that arise in estimating the Fourier transforms  of surface carried measures, we can again reduce to the domain \eqref {restdomain},  where
$
|x_2-cx_1^m|\le \ve x_1^m.
$

\medskip
In a {\it second step}, if $\ka$ denotes again the principal weight associated to the principal edge of $\N(\phi),$ we again perform a dyadic decomposition of the measure $\mu$ over the domain  \eqref {restdomain}, and re-scale in a suitable way. This leads to a phase function 
\begin{equation}\nonumber
\phi(x,\de):= b( \de_1x_1, \de_2x_2)\Big(x_2-x_1^m\om(\de_1x_1)\Big)^2+\de_0 x_1^n \beta(\delta_1x_1),
\end{equation}
where $\de=(\de_0,\,\de_1,\,\de_2)=(2^{-(n\ka_1-1)k},2^{-\ka_1k},2^{-\ka_2k}) $ are small
parameters  which tend to $0$ as $k$ tends to infinity, and where $b( \de_1x_1, \de_2x_2)\sim b(0,0)\ne 0$ and $\beta(0)\ne 0.$

What we then need to prove is the following 
\begin{proposition}\label{s5.1}
Given any point $v=(v_1,v_2)$ such that $v_1\sim 1$ and $v_2=v_1^m\om(0),$ there exists a neighborhood $V$ of $v$ in $(\RR_+)^2$ such that for every cut-off function $\eta\in \D(V),$ the measure 
$\nu_\de$ given by
$$
\laa \nu_\de,f\ra:=\int f(x,\phi (x,\de))\, \eta\x\, dx
$$
satisfies a restriction estimate 
$$
\Big(\int |\widehat f|^2\,d\nu_\de\Big)^{1/2}\le C_{p,\eta} \|f\|_{L^p(\RR^3)},
$$
whenever $p'\ge  2d+2,$ provided $\de$ is sufficiently small.
\end{proposition}
In oder to prove this proposition, we again perform a  dyadic decomposition, this time with respect to the $x_3$-variable. By means of Littlewood-Paley theory, it then turns out that it is sufficient to prove  uniform restriction estimates for the following family of measures
$$
\laa \nu_{\de,j},f\ra:=\int f(x,\phi (x,\de))\, \chi(2^{2j}\phi(x,\de))\eta\x\, dx,
$$
of the form
\begin{equation}\nonumber
\Big(\int |\widehat f|^2\,d\nu_{\de,j}\Big)^{1/2}\le C_{p,\eta} \|f\|_{L^p(\RR^3)}.
\end{equation}
Here, $\chi(x_3) $ is supported on a set where $|x_3|\sim 1.$ If $2^{2j}\de_0\ll 1,$ then it turns out that this localization means in fact again a localization  to a curved rectangle where 
$$|x_1-v_1|< \ve\quad \mbox{ and } \quad |x_2-x_1^m\om(\de_1x_1)|\sim 2^{-j},
$$
but in other cases, it has another meaning.

It turns out that  these estimates require a {\it refined spectral decomposition} of the measures $\nu_{\de,j},$ namely a dyadic decomposition in every dual coordinate $\xi_1,\xi_2,\xi_3.$  Slightly cheating,  for every triple $\La=(\la_1,\la_2,\la_3)$ of dyadic numbers $\la_i=2^{-k_i},$ we therefore define the functions 
$\nu^{\La}_{j}$ by 
\begin{equation}\nonumber
\widehat{\nu^\La_{j}}(\xi)=\chi_1\Big(\frac{\xi_1}{\la_1}\Big)\chi_1\Big(\frac{\xi_2}{\la_2}\Big)\chi_1\Big(\frac{\xi_3}{\la_3}\Big)\widehat{\nu_{\de,j}}(\xi),
\end{equation}
and choose $\chi_1(s)$ supported where $|s|\sim 1$ so that $\nu_{\de,j}=\sum_\La \nu^{\La}_{j},$ where summation is essentially over all these dyadic triples $\La.$ We have here suppressed the dependency of  $\nu^{\La}_{j}$ on the small parameters $\de.$
Note that  $|\xi_i| \sim\la_i $ on the support of $\widehat{\nu^{\la}_{j}}.$
\medskip

For every fixed $\La,$ we then essentially follow again the Stein-Tomas approach, by estimating $\|\widehat{\nu^{\La}_{j}}\|_\infty$ and $\|\nu^{\La}_{j}\|_\infty.$ 

This requires the distinction of various cases, depending on the relative sizes of $\la_1,\la_2$ and $\la_3.$ In the end, it turns out that the most difficult case is where $\la_1\sim\la_2\sim\la_3 $  and $2^{2j}\de_0\sim1,$ and the main result to be proven in this case is the following.

\begin{prop}\label{m2_A}
Let $\phi$ be of type $A_{n-1}$, with $m=2$ and  finite $n\ge 5.$ Then
\begin{equation}\label{restm2}
\sum\limits_{2\le \la_1\sim\la_2\sim\la_3\le 2^{6j}}\, \int_S |\widehat f|^2\,d\nu_j^\La\le C\,2^{\frac 17 j}\, \|f\|^2_{L^{14/11}(\RR^3)},\end{equation}
for all $j\in\NN$ sufficiently big, say $j\ge j_0,$ where the constant $C$ does neither depend on $\de,$ nor on  $j.$ 
\end{prop}

The proof of this result requires yet further refinements. Indeed, the Fourier transform of $\nu^{\La}_{j}$ is an oscillatory integral with complete phase of the form
\begin{eqnarray*}
\Phi(y;\de,j,\xi)&=&\xi_1y_1+\xi_2y_1^2\om(\de_1y_1)+\xi_3\si y_1^n \beta(\de_1 y_1)  \\
&&+2^{-j}\xi_2y_2+\xi_3b^\sharp(y,\de,j)\,y_2^2,\
\end{eqnarray*}
where $\si:=2^{2j}\de_0\sim 1$ and  
 $|b^\sharp(y,\de,j)|\sim 1.$ Notice that if $|\xi_1|\sim|\xi_2|\sim|\xi_2|,$ then $\phi$ may have degenerate critical points, with non-vanishing third derivatives, with respect to the variable $y_1,$  so that we encounter oscillatory integrals of  ``Airy type''.
 
This case requires a {\it further dyadic frequency decomposition with respect  to the distance to certain ``Airy cones.'' } 
For the quite comprehensive details, I refer to \cite{IM-rest1}.

%and in addition somewhat subtle complex interpolation methods in order to capture also the endpoint $p=p_c=14/11.$  For the quite comprehensive details, I refer to \cite{IM-rest2}.
 
 \medskip 
 The case where $2\le d<5$ turns out to be the most delicate one, since on the one hand, we can no longer apply Drury's restriction theorem, which necessitates an even  more refined domain decomposition, and furthermore, somewhat subtle interpolation arguments are needed in many situations in order to handle the endpoint $p=p_c$ (cf. \cite{IM-rest2}).
 
%  \color{red}
%AENDERUNG!!!
%\color{black}


\begin{thebibliography}{9}

\bibitem{arhipov}
 Arhipov, G.\,I.,  Karacuba, A.\,A.,   {\v{C}}ubarikov, V.\,N.,
\newblock Trigonometric integrals.
\newblock {\em Izv. Akad. Nauk SSSR Ser. Mat.}, 43 (1979), 971--1003, 1197 (Russian);
English translation in {\em Math. USSR-Izv.}, 15 (1980), 211--239.

\bibitem{arnold}
 Arnol'd, V.\,I.,
\newblock Remarks on the method of stationary phase and on the {C}oxeter
  numbers.
\newblock {\em Uspekhi Mat. Nauk}, 28 (1973),17--44 (Russsian); English translation 
in {\em Russian Math. Surveys}, 28 (1973), 19--48.


\bibitem{agv}{
Arnol'd, V.I.,  Gusein-Zade, S.M. and  Varchenko, A.N., 
\newblock { Singularities of differentiable maps.  {Vol. II}},  Monodromy and asymptotics of integrals,
 {\em Monographs in Mathematics, 83}.
\newblock Birkh\"auser, Boston Inc., Boston, MA, 1988.}

\bibitem{atiyah}
Atiyah, M., 
\newblock Resolution of singularities and division of distributions.
\newblock {\em Comm. Pure Appl. Math.}, 23 (1970) 145--150.

\bibitem{bak-oberlin-seeger}
Bak, J.G., Oberlin, D.M., Seeger, A.,
\newblock Restriction of Fourier transforms to curves and related oscillatory integrals.
\newblock {\em Amer. J. Math.}, 131 (2) (2009), 277--311.

\bibitem{bernstein-gelfand}
Bernstein, I.N., Gelfand, S.I., 
\newblock Meromorphy of the function $P^\lambda.$
\newblock {\em Funktsional. Anal. i Prilo\v{z}en.}, 3 (1) (1969) 84Ð85.


\bibitem{bierstone-milman1}
Bierstone, E, Milman, P. D., 
\newblock Semianalytic and subanalytic sets.
\newblock {\em Inst. Hautes \'Etudes Sci. Publ. Math.}, 67 (1988), 5--42.
   
   
\bibitem{bierstone-milman2}
Bierstone, E, Milman, P. D., 
\newblock Arc-analytic functions.
\newblock {\em Invent. Math.}, 101 (1990), 411--424.
   


\bibitem{bourgain-circle}
Bourgain, J.,
\newblock Averages in the plane over convex curves and maximal operators.
\newblock {\em J. Analyse Math.}, 47 (1986), 69--85. 

\bibitem{bourgain-cone}
Bourgain, J.,
\newblock Estimates for cone multipliers.
\newblock {\em Geometric aspects of functional analysis (Israel, 1992Ð1994), 41Ð60,}
 \newblock Oper. Theory Adv. Appl. 77 (1995), 41--60. 
   

\bibitem{bourgain-guth}
Bourgain, J., Guth, L., 
\newblock Bounds on oscillatory integral operators.
\newblock {\em C. R. Acad. Sci. Paris, Ser. I}, 349 (2011) 137--141. 

\bibitem{bruna-n-w}
Bruna, J.,  Nagel, A.  and Wainger, S., 
\newblock Convex hypersurfaces and {F}ourier transforms.
\newblock {\em Ann. of Math. (2)}, 127 (1988), no. 2, 333--365.

\bibitem{ccw}
Carbery, C., Christ, M. Wright, J., 
 \newblock  Multidimensional van der Corput and sublevel set estimates.
 \newblock {\em J. Amer. Math. Soc.}, 12 (1999), no. 4, 981--1015.
 
\bibitem{carleson-sjoelin}
Carleson, L.,   P. Sj\"olin, P., 
\newblock Oscillatory integrals and a multiplier problem for the disc.
\newblock {\em Studia Math.}, 44 (1972), 287--299.
 
\bibitem{cgp}
Collins, T.,  Greenleaf, A., Pramanik, M.,   
 \newblock A multi-dimensional resolution of singularities with applications to analysis.
 \newblock {\em Preprint, arXiv:1007.0519.}
 
 \bibitem{vdC}
van der Corput, J.\,G., 
\newblock Zahlentheoretische Absch\"atzungen.
\newblock {\em Math. Ann.}, 84 (1921),  53--79.

 
 \bibitem{cowling-d-m}
Cowling, M., Disney, S., Mauceri, G., M\"uller, D.,
\newblock Damping oscillatory integrals.
\newblock {\em Invent. Math.}, 101 (1990), 237--260.

 
 \bibitem{cowling-mauceri2}
Cowling, M. and Mauceri, G.,
\newblock Inequalities for some maximal functions. {II}.
\newblock {\em Trans. Amer. Math. Soc.}, 296 (1986), no. 1, 341--365.

\bibitem{cowling-mauceri1}
Cowling, M. and Mauceri, G.,
\newblock Oscillatory integrals and {F}ourier transforms of surface carried
  measures.
\newblock {\em Trans. Amer. Math. Soc.}, 304 (1987), no.1, 53--68.

 \bibitem{dendrinos-m}
Dendrinos, S., M\"uller, D., 
\newblock Uniform estimates for the local restriction of the Fourier transform to curve.
\newblock {\em Trans. Amer. Math. Soc.}, to appear.


 
 
\bibitem{dns}
Denef, J., Nicaise, J.,  Sargos, P.,
\newblock  Oscillatory integrals and Newton polyhedra.
\newblock {\em J. Anal. Math.}, 95 (2005), 147--172.

\bibitem{domar}
Domar,Y., 
\newblock On the {B}anach algebra {$A(G)$} for smooth sets {$\Gamma \subset
  \RR^n$}.
\newblock {\em Comment. Math. Helv.}, 52, no. 3 (1977), 357--371.


\bibitem{drury}
S. W. Drury,
\newblock Restrictions of Fourier transforms to curves,
\newblock {\em Ann. Inst. Fourier.}, (35) (1985), 117--123.



\bibitem{duistermaat}
 Duistermaat, J.\,J.,
\newblock Oscillatory integrals, {L}agrange immersions and unfolding of
  singularities.
\newblock {\em Comm. Pure Appl. Math.}, 27 (1974), 207--281.

\bibitem{fefferman}
Fefferman, C.,
\newblock Inequalities for strongly singular convolution operators.
\newblock {\em Acta Math.}, (1970), 9--36.

\bibitem{grafakos}
Grafakos, L.,
\newblock {Modern Fourier analysis}.
 {\em Graduate Texts in Mathematics} 250.
\newblock Springer, New York, 2009.


\bibitem{greenblatt3}
Greenblatt, M., 
\newblock Newton polygons and local integrability of negative powers of
   smooth functions in the plane.
\newblock {\em Trans. Amer. Math. Soc.}, 358 (2006), no. 2, 657--670.


\bibitem{greenblatt-resol1}
Greenblatt, M., 
\newblock An elementary coordinate-dependent local resolution of singularities and applications.
\newblock {\em  J. Funct. Anal.},  255 (2008), no. 8, 1957--1994.

\bibitem{greenblatt4}
Greenblatt, M., 
\newblock The asymptotic behavior of degenerate oscillatory integrals in two dimensions.
\newblock {\em J. Funct. Anal.}, 257 (2009), no. 6, 1759--1798.

\bibitem{greenblatt-resol2}
Greenblatt, M., 
\newblock Resolution of singularities, asymptotic expansions of integrals and related phenomena.
\newblock {\em   J. Anal. Math.},   111 (2010), 221--245.

\bibitem{greenblatt1}
Greenblatt, M., 
\newblock Oscillatory integral decay, sublevel set growth, and the Newton polyhedron.
\newblock {\em Math. Ann.}, 346 (2010),  857--895.


\bibitem{greenblatt-max1}
Greenblatt, M., 
\newblock $L^p$ boundedness of maximal averages over hypersurfaces in $\RR^3.$ 
\newblock {\em Preprint}

\bibitem{greenblatt-max2}
Greenblatt, M., 
\newblock Maximal averages over hypersurfaces and the Newton polyhedron.
\newblock {\em Preprint}, arXiv:1002.0109.




\bibitem{greenleaf}
Greenleaf, A.,
\newblock Principal curvature and harmonic analysis.
\newblock {\em Indiana Univ. Math. J.}, 30(4) (1981), 519--537.

\bibitem{hironaka} Hironaka, H.,
\newblock Resolution of singularities of an algebraic variety over a field of characteristic zero I, II.
\newblock {\em  Ann. Math.,} (2), 79 (1964), 109--326.

%\bibitem{hormander-book} L. H\"ormander
%\newblock The analysis of linear partial differential operators. I. Distribution theory and Fourier analysis. Reprint of the second (1990) edition [Springer, Berlin; MR1065993 (91m:35001a)]. 
%\newblock {\em Classics in Mathematics. Springer-Verlag, Berlin,} 2003. x+440 pp. ISBN: 3-540-00662-1 35-02

\bibitem{hoermander}
H\"ormander, L.,
\newblock Oscillatory integrals and multipliers on FLp.
\newblock {\em Ark. Mat.}, 11 (1973), 1--11.

\bibitem{ikm}
Ikromov, I.A., Kempe, M. and   M{\"u}ller, D.,
\newblock Damped oscillatory integrals and boundedness of maximal operators
  associated to mixed homogeneous hypersurfaces.
\newblock {\em Duke Math. J.}, 126 (2005), no. 3, 471--490.


\bibitem{IKM-max} Ikromov, I.\,A., Kempe, M.,  M\"uller, D.,
\newblock Estimates for maximal functions  associated to hypersurfaces in $\bR^3$  and related problems of harmonic analysis.  
\newblock {\em Acta Math.} 204 (2010), 151--271.

\bibitem{IM-ada}
Ikromov, I.\,A.,  M\"uller, D., 
\newblock On adapted coordinate systems.
\newblock {\em   Trans. Amer. Math. Soc.,} 363 (2011), no. 6, 2821--2848. 


\bibitem{IM-uniform}
Ikromov, I.\,A.,  M\"uller, D., 
\newblock Uniform estimates for the Fourier transform of surface carried measures in  $\bR^3$ and an application to Fourier restriction.
\newblock {\em    J. Fourier Anal. Appl.,} 17 (2011), no. 6, 1292--1332.


\bibitem{IM-rest1}
Ikromov, I.\,A.,  M\"uller, D., 
\newblock $L^p$-$L^2$  Fourier restriction  for  hypersurfaces in  $\bR^3:$ Part I.
\emph{Preprint; ArXiv:  http://arxiv.org/abs/1208.6090}

\bibitem{IM-rest2}
Ikromov, I.\,A.,  M\"uller, D., 
\newblock $L^p$-$L^2$  Fourier restriction  for  hypersurfaces in  $\bR^3:$ Part II.
\emph{In preparation.}

%\bibitem{iosevich}
%Iosevich, A., 
%Fourier transform, $L^2$ restriction theorem, and scaling.
%\emph{Boll. Unione Mat. Ital. Sez. B Artic. Ric. Mat.} (8) 2 (1999), 383--387.

\bibitem{iosevich-curves}
Iosevich, A., 
\newblock Maximal operators associated to families of flat curves in the plane.
\newblock {\em Duke Math. J.}, 76 (1994), no. 2, 633--644.

\bibitem{iosevich-sawyer1}
Iosevich, A. and Sawyer, E.,
\newblock Oscillatory integrals and maximal averages over homogeneous surfaces.
\newblock {\em Duke Math. J.}, 82 (1996), no. 1, 103--141.

\bibitem{iosevich-sawyer2}
Iosevich, A. and Sawyer, E.,
\newblock Maximal averages over surfaces.
\newblock {\em Adv. Math.}, 132 (1997), no. 1, 46--119.

\bibitem{io-sa-seeger}
Iosevich, A.,  Sawyer, E. and Seeger, A.,
\newblock On averaging operators associated with convex hypersurfaces of finite
  type.
\newblock {\em J. Anal. Math.}, 79 (1999), 159--187.

\bibitem{karpushkin}
 Karpushkin, V.\,N.,
\newblock A theorem on uniform estimates for oscillatory integrals with a phase
  depending on two variables. 
  \newblock {\em Trudy Sem. Petrovsk.} 10 (1984), 150--169, 238 (Russian);
  English translation in 
{\em  J. Soviet Math.}, 35 (1986), 2809--2826.

\bibitem{magyar}
Magyar, A., 
\newblock On Fourier restriction and the Newton polygon.
\newblock {\em  Proceedings Amer. Math. Soc.} 137 (2009), 615--625.

\bibitem{mockenhaupt}
Mockenhaupt, G., 
\newblock A note on the cone multiplier.
\newblock {\em Proc. Amer. Math. Soc.}, 117 (1993), no. 1, 145--152. 

\bibitem{mockenhaupt-seeger-sogge2}
Mockenhaupt, G.,  Seeger, A.  and Sogge, C.D.,
\newblock Wave front sets, local smoothing and {B}ourgain's circular maximal
  theorem.
\newblock {\em Ann. of Math. (2)}, 136 (1992), no. 1, 207--218.


\bibitem{mockenhaupt-seeger-sogge}
Mockenhaupt, G.,  Seeger, A.  and Sogge, C.D.,
\newblock Local smoothing of {F}ourier integral operators and
  {C}arleson-{S}j\"olin estimates.
\newblock {\em J. Amer. Math. Soc.}, 6 (1993), no. 1, 65--130.



\bibitem{nagel-seeger-wainger}
Nagel, A.,  Seeger, A.  and Wainger, S.,
\newblock Averages over convex hypersurfaces.
\newblock {\em Amer. J. Math.}, 115 (1993), no. 4, 903--927.

\bibitem{parusinski1}
Parusi\'nski, A.,
\newblock Subanalytic functions.
\newblock {\em Trans. Amer. Math. Soc.}, 344 (1994), 583--595.

\bibitem{parusinski2}
Parusi\'nski, A.,
\newblock On the preparation theorem for subanalytic functions.
\newblock {\em  New developments
in singularity theory (Cambridge 2000), Nato Sci. Ser. II Math. Phys. Chem., 21}, 
Kluwer Academic Publisher, Dordrecht (2001), 193--215.



\bibitem{phong-stein}
Phong, D.H. and  Stein, E.M., 
\newblock The {N}ewton polyhedron and oscillatory integral operators.
\newblock {\em Acta Math.}, 179 (1997), no. 1, 105--152.


\bibitem{PSS} Phong, D.\,H., Stein, E.\,M.,  Sturm, J.\,A.,
\newblock On the growth and stability of real-analytic functions.
\newblock {\em Amer. J. Math.}, 121 (1999), no. 3, 519-554.

\bibitem{randol}
Randol, B., 
\newblock On the asymptotic behavior of the {F}ourier transform of the
  indicator function of a convex set.
\newblock {\em Trans. Amer. Math. Soc.}, 139 (1969), 279--285.

\bibitem{schulz-un}
Schulz, H., 
\newblock On the decay of the Fourier transform of measures on hypersurfaces, generated by radial functions, and related restriction theorems.
 \newblock {\em unpublished preprint, 1990.} 
 
\bibitem{schulz}
Schulz, H., 
\newblock Convex hypersurfaces of finite type and the asymptotics of their
  {F}ourier transforms.
\newblock {\em Indiana Univ. Math. J.}, 40 (1991), 1267--1275.



\bibitem{sogge-stein}
Sogge, C.D.  and Stein, E.M., 
\newblock Averages of functions over hypersurfaces in {$\RR^n$}.
\newblock {\em Invent. Math.}, 82 (1985), no. 3, 543--556.

\bibitem{stein-sphere}
Stein, E.M.,
\newblock Maximal functions. {I}. {S}pherical means.
\newblock {\em Proc. Nat. Acad. Sci. U.S.A.}, 73 (1976), no. 7, 2174--2175.

\bibitem{stein-book}
Stein, E.M.,
\newblock {Harmonic analysis: Real-variable methods, orthogonality, and
  oscillatory integrals}. {\em Princeton Mathematical Series} 43.
\newblock Princeton University Press, Princeton, NJ, 1993.

%\bibitem{stein-weiss}
%Stein, E.\,M., Weiss, G,,
%\newblock {Introduction to Fourier analysis on Euclidean spaces. }
%{\em Princeton Mathematical Series} No. 32. 
%\newblock Princeton University Press, Princeton, N.J., 1971.






\bibitem{strichartz} Strichartz, R.\,S.,  Restrictions of Fourier transforms to quadratic surfaces and decay of solutions of wave equations. 
\newblock {\em Duke Math. J.}, 44  (1977), 705--714.

\bibitem{sussman}  Sussmann, H.J,.   Real-analytic desingularization and subanalytic sets: an elementary
approach.
\newblock {\em Trans. Amer. Math. Soc.}, 317 (1990), 417--461.

\bibitem{svensson}
Svensson, I.,
\newblock Estimates for the {F}ourier transform of the characteristic function
  of a convex set.
\newblock {\em Ark. Mat.}, 9 (1971), 11--22.



\bibitem{Va}
Varchenko, A.\,N.,
\newblock Newton polyhedra and estimates of oscillating integrals.
\newblock {\em Funkcional. Anal. i Prilo\v zen}, 10 (1976), 13--38 (Russian); English translation in 
\newblock {\em Funktional Anal. Appl.}, 18 (1976), 175--196.

\bibitem{zygmund}
Zygmund, A.,
\newblock On Fourier coefficients and transforms of functions of two variables.
\newblock {\em Studia Math.}, 9 (1971), (1974), 189--201.



\end{thebibliography}
\end{document}